\documentclass[12pt]{article}
\usepackage{graphicx} % Required for inserting images
\usepackage{hyperref}
\usepackage{amsmath,amscd,mathtools,cases,amsthm,amsopn,amssymb,soul,cleveref,mathabx,relsize,tikz,cleveref}
\usepackage{longtable,geometry}

\usepackage{thmtools}
\usepackage{thm-restate}

\newcommand\Tstrut{\rule{0pt}{12ex}}       % "top" strut
\newcommand\Bstrut{\rule[-11.5ex]{0pt}{0pt}} % "bottom" strut
\newcommand{\TBstrut}{\Tstrut\Bstrut} % top&bottom struts

\newtheorem{theorem}[subsection]{Theorem}
\newtheorem{corollary}[subsection]{Corollary}
\newtheorem{observation}[subsection]{Observation}
\newtheorem{proposition}[subsection]{Proposition}
\newtheorem{lemma}[subsection]{Lemma}

\theoremstyle{definition}
\newtheorem{definition}[subsection]{Definition}
\theoremstyle{definition}

\newtheorem{example}[subsection]{Example}

\DeclareMathOperator{\mr}{mr}
\DeclareMathOperator{\rank}{rank}
\DeclareMathOperator{\nullity}{null}
\DeclareMathOperator{\N}{\mathcal{N}}

\def\GraphSix#1"{\begingroup(\itshape graph6:\nolinebreak[3]\ \verb"\aftergroup\endgroup\aftergroup)}

\crefname{figure}{figure}{}
\Crefname{figure}{Figure}{}

\title{The Classification of Graphs on $8$ vertices with Coinciding Zero Forcing number and Maximum Nullity}

\begin{document}

\author{Wayne Barrett \thanks{Brigham Young University, 155 East 1230 North, Provo, UT 84602, wb@mathematics.byu.edu} \and Mark Hunnell \thanks{Winston-Salem State University, 601 S Martin Luther King Jr Dr, Winston-Salem NC 27110, hunnellm@wssu.edu} \and John Hutchens \thanks{University of San Francisco, 2130 Fulton St, San Francisco CA 94117, johnd.hutchens@gmail.com} \and John Sinkovic \thanks{Brigham Young University-Idaho, 101 E Viking St, Rexburg, ID 83460, sinkovicj@byui.edu} }
\date{}

\maketitle
\begin{abstract}
    We study the minimum rank of a (simple, undirected) graph, which is the minimum rank among all matrices in a space  determined by the graph. We determine the exact set of graphs on eight vertices for which the nullity of a minimum rank matrix does not coincide with a bound determined by the zero forcing number of a graph. Although our goal was to determine which eight-vertex graphs satisfy maximum nullity equal to the zero forcing number, we also established several additional methods to assist in the computation of minimum rank for general graphs. 
\end{abstract}



\section{Introduction}

Determining the maximum multiplicity of an eigenvalue among all symmetric matrices in a space of matrices defined by a graph (definitions to follow) is an important relaxation of the Inverse Eigenvalue Problem for a Graph (IEPG). 
 In this article we describe new approaches for this problem on small graphs, and completely determine the minimum rank of all graphs with eight or fewer vertices.  Work has continued on the minimum rank problem for a graph since \cite{deloss2008table}, though some of the techniques established since then have not been implemented in code; we provide the relevant update to the algorithm described therein.

Given a simple undirected graph on $n$ vertices, let $S(G)$ be the set of all $n \times n$ matrices $A$ such that $a_{ij} \ne 0$ if and only if $ij \in E(G)$.  There is no condition on the diagonal entries of $A$ (which allows for diagonal translation).  Let 

\[ \mr(G) = \min \{ \ \rank A \ | \ A \in S(G) \ \},  \] 

and let 

\[ M(G)=\max \{ \ \nullity A \  | \ A \in S(G) \ \}.  \] 

It is immediate that $\mr(G) + M(G) =n$, and thus for fixed $n$ the problems of determining $\mr(G)$ and $M(G)$ are equivalent.  A matrix $A \in S(G)$ satisfying $\rank A = \mr(G)$ is called a witness of $\mr(G)$.  Determining the minimum rank /maximum nullity of a graph has been studied extensively since \cite{nylen1996minimum}.  For most graphs $\mr(G) \ (\text{equivalently, }M(G))$ is difficult to compute.  

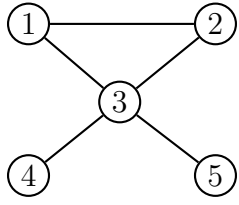
\begin{figure}
    \centering
\begin{tikzpicture}[scale=1.5,thick]
\tikzstyle{every node}=[minimum width=0pt, inner sep=2pt, circle]
    \draw (-1.33,2.07) node[draw] (1) {$1$};
    \draw (0.31,2.07) node[draw] (2) {$2$};
    \draw (-0.53,1.38) node[draw] (3) {$3$};
    \draw (-1.33,0.73) node[draw] (4) {$4$};
    \draw (0.31,0.73) node[draw] (5) {$5$};
    \draw  (1) edge (3);
    \draw  (2) edge (3);
    \draw  (3) edge (5);
    \draw  (3) edge (4);
    \draw  (1) edge (2);
\end{tikzpicture}
    \caption{The campstool graph $W$.}
    \label{fig:campstool}
\end{figure}

In \cite{aim2008zero} a close connection was found between $M(G)$ and a graph parameter called the zero forcing number.  This is defined as follows.  Each vertex of a graph is colored blue or white.  We then define the standard color change rule:
\begin{itemize}
    \item[] \textbf{Standard color change rule}: If a blue vertex $b$ has a unique white neighbor $w$, then change the color of $w$ to blue.
\end{itemize}
If $w$ is the unique white neighbor of a blue vertex $b$, we say $b$ forces $w$ and denote this $b \rightarrow w$. A subset $Z \subset V(G)$ is called a zero forcing set if when the vertices of $Z$ are colored blue, all vertices in $G$  can be made blue by successively applying the standard color change rule.  The zero forcing number $Z(G)$ of a graph is the minimum size of a zero forcing set.  For example, in the campstool graph in \Cref{fig:campstool}, no single vertex is a zero forcing set, but $\{1,4\}$ is a zero forcing set.  Therefore $Z({\rm campstool})=2$.  For most graphs $G$, calculating $Z(G)$ is NP-hard, although several known algorithms exist.  In general, at various points in the application of the color change rule it may be possible for more than one blue vertex to force a white vertex.  Given a zero forcing set $B$, we can choose a chronological set of forces $F$ such that each force in the set is valid given the preceding forces.  A chronological set of forces $F$ induces a chain $C_v$ for each $v \in B$, where $C_v = \{ v=v_0, v_1, \dots, v_k \}$ and $v_{i-1} \rightarrow v_i \in F$ for each $i \in \{ 1,2, \dots k\}$.  We allow chains to consist of one vertex when $v\in B$ does not perform a force in $F$.  The collection of chains $\{ C_v \}_{v\in B}$ is called a set of zero forcing chains of $B$.

An exact algorithm for calculating the minimum rank of a graph is known, but the number of variables makes the computation cost prohibitive. In \cite{aim2008zero}, a purely combinatorial bound on the maximum nullity was proven:

\begin{theorem}
     For any graph $G$, $M(G) \le Z(G)$.
\end{theorem}

A natural question is when equality holds, and a satisfactory answer has remained elusive.  Early on in the study of zero forcing, it was also shown in \cite{aim2008zero} that:

\begin{theorem}\label{nleq6}
    If $G$ is a graph on 6 or fewer vertices, then $M(G)=Z(G)$.
\end{theorem}
\medskip

These results generated a lot of interest in $Z(G)$.  They were extended in \cite{deloss2008table}, where the minimum rank of all 1044 graphs on 7 vertices was calculated. By examination of their table it can be seen that $M(G)=Z(G)$ for all graphs of order $n\leq 7$.  In addition, there are several known graph families for which $M(G) =Z(G)$, such as the complete graphs $K_n$, complete bipartite graphs, cycles, and trees \cite{aim2008zero}.

The equality $M(G)=Z(G)$ no longer holds for eight-vertex graphs, and one principal goal of this paper is to identify  the connected eight-vertex graphs for which $M(G) < Z(G)$.

We are now ready to state our main result.

\begin{restatable}{theorem}{mainthm}\label{mainthm}  Let $\mathcal{E} = \{ E_1, E_2, \dots , E_7\}$ be the set of graphs given in \Cref{fig:E1E7}. A connected graph $G$ on 8 vertices satisfies $M(G) < Z(G)$ if and only if $G \in \mathcal{E}$.
\end{restatable}

\begin{figure}
    \centering

\begin{tikzpicture}[scale=1.0,thick]
    \tikzstyle{every node}=[minimum width=0pt, inner sep=2pt, circle]

    \node at (-.25,2) (E2) {$E_2$};
	\draw (1.5,-0.61) node[draw] (0) {};
	\draw (0.05,0.82) node[draw] (1) {};
	\draw (1.45,2.23) node[draw] (2) {};
	\draw (0.93,1.25) node[draw] (3) {};
	\draw (1.82,0.8) node[draw] (4) {};
	\draw (1.99,1.67) node[draw] (5) {};
	\draw (0.92,0.38) node[draw] (6) {};
	\draw (0.45,1.83) node[draw] (7) {};
	\draw  (4) edge (6);
	\draw  (0) edge (6);
	\draw  (1) edge (6);
	\draw  (0) edge (2);
	\draw  (2) edge (3);
	\draw  (1) edge (3);
	\draw  (3) edge (4);
	\draw  (0) edge (4);
	\draw  (0) edge (1);
	\draw  (1) edge (7);
	\draw  (2) edge (7);
	\draw  (2) edge (5);
	\draw  (4) edge (5);
\begin{scope}[shift={(-1,-1.3)},scale=1.5]
    \draw (-0.59,2.35) node[draw] (0) {  };
    \draw (-0.59,1.69) node[draw] (1) { };
    \draw (-1.29,1.69) node[draw] (2) {};
    \draw (0.09,1.69) node[draw] (3) {};
    \draw (-0.92,1.06) node[draw] (4) {};
    \draw (-0.27,1.06) node[draw] (5) {};
    \draw (-1.29,0.49) node[draw] (6) {};
    \draw (0.09,0.49) node[draw] (7) {};
    \draw  (0) edge (1);
    \draw  (1) edge (3);
    \draw  (1) edge (2);
    \draw  (1) edge (4);
    \draw  (1) edge (5);
    \draw  (3) edge (5);
    \draw  (4) edge (5);
    \draw  (2) edge (4);
    \draw  (4) edge (6);
    \draw  (5) edge (7);
    \node (a) at (-1,2.2) {$E_1$};
    \end{scope}
% \end{tikzpicture}

% % Graph6: GzcmR?
% $E_3$ \begin{tikzpicture}[scale=1.5,thick]
%     \tikzstyle{every node}=[minimum width=0pt, inner sep=2pt, circle]
\begin{scope}[shift={(3,0)}]
        \node at (-.25,2) (E2) {$E_3$};

	\draw (1.5,-0.61) node[draw] (0) {};
	\draw (0.05,0.82) node[draw] (1) {};
	\draw (1.45,2.23) node[draw] (2) {};
	\draw (0.93,1.25) node[draw] (3) {};
	\draw (1.82,0.8) node[draw] (4) {};
	\draw (1.99,1.67) node[draw] (5) {};
	\draw (0.92,0.38) node[draw] (6) {};
	\draw (0.45,1.83) node[draw] (7) {};
	\draw  (2) edge (5);
	\draw  (4) edge (5);
	\draw  (2) edge (3);
	\draw  (2) edge (7);
	\draw  (1) edge (2);
	\draw  (0) edge (2);
	\draw  (1) edge (7);
	\draw  (1) edge (3);
	\draw  (3) edge (4);
	\draw  (0) edge (4);
	\draw  (4) edge (6);
	\draw  (0) edge (6);
	\draw  (1) edge (6);
	\draw  (0) edge (1);

 \end{scope}
% \end{tikzpicture}
% %graph6: GrsmR?
% \hskip 1in $E_4$ \begin{tikzpicture}[scale=1.5,thick]
%     \tikzstyle{every node}=[minimum width=0pt, inner sep=2pt, circle]
\begin{scope}[shift={(6,0)}]
        \node at (-.25,2) (E2) {$E_4$};
	\draw (1.5,-0.61) node[draw] (0) {};
	\draw (0.05,0.82) node[draw] (1) {};
	\draw (1.45,2.23) node[draw] (2) {};
	\draw (0.93,1.25) node[draw] (3) {};
	\draw (1.82,0.8) node[draw] (4) {};
	\draw (1.99,1.67) node[draw] (5) {};
	\draw (0.92,0.38) node[draw] (6) {};
	\draw (0.45,1.83) node[draw] (7) {};
	\draw  (0) edge (1);
	\draw  (0) edge (6);
	\draw  (0) edge (4);
	\draw  (0) edge (2);
	\draw  (1) edge (6);
	\draw  (4) edge (6);
	\draw  (1) edge (4);
	\draw  (1) edge (7);
	\draw  (2) edge (7);
	\draw  (2) edge (5);
	\draw  (4) edge (5);
	\draw  (2) edge (3);
	\draw  (3) edge (4);
	\draw  (1) edge (3);
 \end{scope}
% \end{tikzpicture}

% \vskip 1cm

% %graph6: GzkmR?
% \qquad $E_5$\begin{tikzpicture}[scale=1.5,thick]
%     \tikzstyle{every node}=[minimum width=0pt, inner sep=2pt, circle]

\begin{scope}[shift={(-1,-4)}]
        \node at (-.25,2) (E2) {$E_5$};
	\draw (1.5,-0.61) node[draw] (0) {};
	\draw (0.05,0.82) node[draw] (1) {};
	\draw (1.45,2.23) node[draw] (2) {};
	\draw (0.93,1.25) node[draw] (3) {};
	\draw (1.82,0.8) node[draw] (4) {};
	\draw (1.99,1.67) node[draw] (5) {};
	\draw (0.92,0.38) node[draw] (6) {};
	\draw (0.45,1.83) node[draw] (7) {};
	\draw  (4) edge (5);
	\draw  (2) edge (5);
	\draw  (3) edge (4);
	\draw  (0) edge (4);
	\draw  (0) edge (6);
	\draw  (4) edge (6);
	\draw  (0) edge (2);
	\draw  (1) edge (2);
	\draw  (1) edge (3);
	\draw  (2) edge (3);
	\draw  (2) edge (7);
	\draw  (1) edge (7);
	\draw  (0) edge (1);
	\draw  (1) edge (6);
	\draw  (2) edge (4);
  \end{scope}
% \end{tikzpicture}
% %graph6: GzsmR?
% \hskip 1in $E_6$
% \begin{tikzpicture}[scale=1.5,thick]
%    \tikzstyle{every node}=[minimum width=0pt, inner sep=2pt, circle]

\begin{scope}[shift={(2,-4)}]
        \node at (-.25,2) (E2) {$E_6$};
	\draw (1.5,-0.61) node[draw] (0) {};
	\draw (0.05,0.82) node[draw] (1) {};
	\draw (1.45,2.23) node[draw] (2) {};
	\draw (0.93,1.25) node[draw] (3) {};
	\draw (1.82,0.8) node[draw] (4) {};
	\draw (1.99,1.67) node[draw] (5) {};
	\draw (0.92,0.38) node[draw] (6) {};
	\draw (0.45,1.83) node[draw] (7) {};
	\draw  (0) edge (6);
	\draw  (0) edge (4);
	\draw  (0) edge (1);
	\draw  (1) edge (6);
	\draw  (0) edge (2);
	\draw  (1) edge (2);
	\draw  (4) edge (6);
	\draw  (4) edge (5);
	\draw  (2) edge (5);
	\draw  (2) edge (3);
	\draw  (3) edge (4);
	\draw  (1) edge (3);
	\draw  (1) edge (7);
	\draw  (2) edge (7);
	\draw  (1) edge (4);
  \end{scope}
% \end{tikzpicture}

% \vskip 1cm
% %graph6: Gz{mR
% \quad $E7$
% \begin{tikzpicture}[scale=1.5,thick]
%     \tikzstyle{every node}=[minimum width=0pt, inner sep=2pt, circle]
\begin{scope}[shift={(5,-4)}]
        \node at (-.25,2) (E2) {$E_7$};
	\draw (1.5,-0.61) node[draw] (0) {};
	\draw (0.05,0.82) node[draw] (1) {};
	\draw (1.45,2.23) node[draw] (2) {};
	\draw (0.93,1.25) node[draw] (3) {};
	\draw (1.82,0.8) node[draw] (4) {};
	\draw (1.99,1.67) node[draw] (5) {};
	\draw (0.92,0.38) node[draw] (6) {};
	\draw (0.45,1.83) node[draw] (7) {};
	\draw  (0) edge (6);
	\draw  (0) edge (1);
	\draw  (1) edge (4);
	\draw  (0) edge (4);
	\draw  (4) edge (6);
	\draw  (4) edge (5);
	\draw  (2) edge (5);
	\draw  (2) edge (4);
	\draw  (1) edge (2);
	\draw  (3) edge (4);
	\draw  (2) edge (3);
	\draw  (1) edge (3);
	\draw  (2) edge (7);
	\draw  (1) edge (7);
	\draw  (0) edge (2);
	\draw  (1) edge (6);
 \end{scope}
\end{tikzpicture}
    \caption{The set of graphs $\mathcal{E}$ in \Cref{mainthm}.}
    \label{fig:E1E7}
\end{figure}
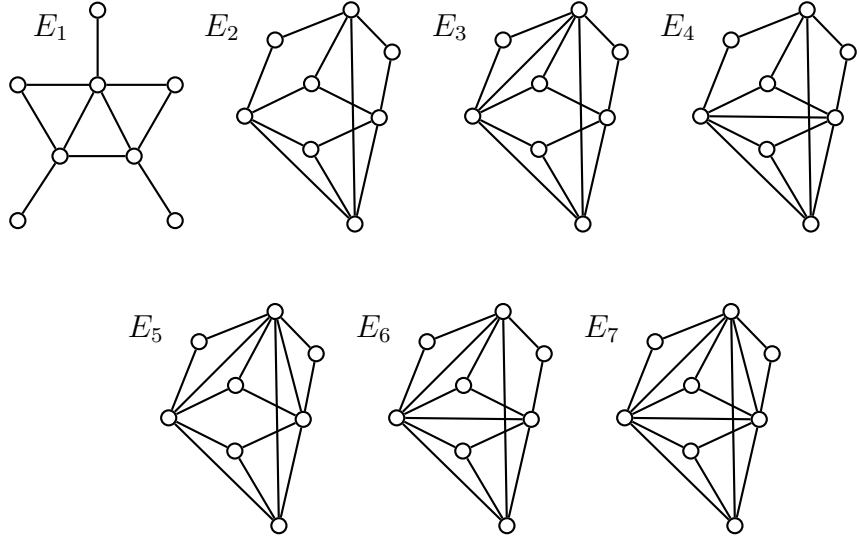

 We now proceed to describe the organization of this article. In \Cref{prelims}, we describe some of the basic results for studying the minimum rank of a graph. \Cref{kpathssec} establishes a general result that simplifies the computation of the minimum rank of some graphs, particularly a subset of graphs on eight vertices for which previously implemented bounds were inconclusive. \Cref{x7vertices} simplifies the proof that all graphs of order $n\leq 7$ satisfy $M(G) = Z(G)$. In \Cref{CutVertex2Sep} we discuss techniques for computing $M(G)$ for graphs with a cut-vertex or a $2$-separation.  \Cref{liftingsec} describes an implementation of using a known witness for $\mr(G-v)$ to determine a witness of $\mr(G)$ under suitable conditions.  A detailed description of the algorithm used to compute $M(G)$ for all but 31 of the connected graphs on $8$ vertices is described in \Cref{algo-details}. Finally, \Cref{final31} describes the technique used to compute a witness for $\mr(G)$ for graphs $G$ that were not amenable to other methods.
 
 \section{Preliminaries}
 \label{prelims}

 Unless otherwise stated, all graphs are simple and undirected. We denote the set of vertices of a graph $G$ by $V(G)$ and its edges by $E(G)$. We denote the complete graph on $n$ vertices by $K_n$.

 A vertex $v \in V(G)$ is a cut vertex if $G-v$ has more connected components than $G$.  If $G$ has a cut vertex $v$, it is natural to consider $G$ as composed of two graphs $G_1, G_2$ where one vertex of each is identified and denoted by $v$.  We denote this decomposition by $G=G_1 \oplus_v G_2$.

 We recall some additional parameters that we will need.  Given a connected graph $G\ne K_n$ on $n$ vertices, its vertex connectivity $\kappa(G)$ is the minimum number of vertices in $V(G)$ whose removal disconnects $G$, while $\kappa(K_n)=n-1$.  A real symmetric $n \times n$ matrix $A$ satisfies the Strong Arnold Hypothesis provided there is no nonzero symmetric $n \times n$ matrix $X$ such that $A X = 0$, $A\circ X=0$, and $I\circ X=0$ where $\circ$ denotes the entry-wise product. The parameter $\xi(G)$ is defined to be the maximum nullity among all real symmetric matrices such that $A \in S(G)$ and $A$ satisfies the Strong Arnold Hypothesis.  The Colin de Verdiere parameter $\mu(G)$ is the largest nullity of any real symmetric matrix $M \in S(G)$ such that:
\begin{enumerate}
\item $M$ satisfies the Strong Arnold Hypothesis,
\item $M$ has exactly one negative eigenvalue, of multiplicity 1, 
\item for all $i, j$ with $i \ne j,  m_{ij} < 0$ if $i$ and $j$ are adjacent   
and $m_{ij} = 0$ if $i$ and $j$ are nonadjacent.
\end{enumerate}

We will make use of the following results.
\begin{theorem} \cite{barioli2013parameters}\label{kappa_xi_M_Z}
For any connected graph $G$,
\begin{equation}
    \kappa(G) \le \xi(G) \le M(G) \le Z(G). 
\end{equation}
\end{theorem}

We will elaborate on the properties of $\mu(G)$ in \Cref{cdv}.

 If a graph $G$ has a cut vertex $v$, then its minimum rank is expressible in terms of its induced subgraphs.  Since we are interested in graphs with eight vertices, the minimum rank of all induced subgraphs is known.
\begin{theorem}
\cite{hsieh2001minimum}\cite{barioli2004computation}  Let $G$ be a graph with a cut vertex $v$ and suppose that $G=G_1\oplus_v G_2$.  Then

\begin{equation}
\mr(G) = \min\{\mr(G_1)+\mr(G_2), \mr(G_1-v)+\mr(G_2-v)+2\}. 
\end{equation}
\end{theorem}

Rewriting this identity in terms of $M$ gives the identity 

\begin{equation}\label{cutvertM}
M(G)=\max \{M(G_1) + M(G_2), M(G_1-v)+M(G_2-v)\}-1.
\end{equation}

Thus, for any graph with a cut vertex, $M(G)$ can be determined from $M(H)$ for four proper induced subgraphs of $G$. A generalization of this formula exists, but requires the following definition.\medskip
 
\begin{definition}
    A 2-separation of a graph $G=(V,E)$ is a pair of subgraphs $G_1=(V_1,E_1)$ and $G_2=(V_2,E_2)$ such that $V_1 \cup V_2\ = \ V, \ E_1 \cup E_2 \ =\ E, \  |V_1 \cap V_2|=2,$ and $E_1 \cap E_2 = \emptyset$.  
\end{definition} 

If a graph $G$ has a 2-separation, then its maximum nullity is determined exactly by a formula.

\begin{theorem} \cite{van2008maximum}
    Let $G$ be a graph with a $2$-separation $(G_1,G_2)$ with $R=\{r_1,r_2\} = V(G_1) \cap V(G_2)$.  Let $H_1$ and $H_2$ be obtained from $G_1$ and $G_2$, respectively, by inserting an edge between $r_1$ and $r_2$, and let $\overline{G_i}$, i=1,2, be the  graph obtained from $G_i$ by identifying $r_1$ and $r_2$.  Then
\setcounter{equation}{3}    
\begin{equation}\label{2sepM2}
\begin{split}
    M(G) = \max \{ &M(G_1) + M(G_2) -2, M(H_1) + M(H_2) - 2, \\
    &M(\overline{G_1}) + M(\overline{G_2}) - 2, M(G_1 - r_1) + M(G_2-r_1) - 2, \\
    &M(G_1-r_2) + M(G_2 - r_2) - 2, M(G_1 - R) + M(G_2-R)-2 \}.
    \end{split}
\end{equation}
\end{theorem}

Similar to a graph with a cut vertex, $M(G)$ for a graph $G$ with a 2-separation can be determined from $M(H)$ for several related graphs $H$ associated with $G$.  

Another tool we employ involves identifying dominating vertices whose deletion preserves equality between the maximum nullity and zero forcing number.

\begin{proposition}\label{domvertexprop}
Suppose that a graph $G$ has a dominating vertex $v$.  Then \\$M(G-v)=Z(G-v)$ implies $M(G)=Z(G)$.
\end{proposition}
\begin{proof}
By Corollary 7.2 in  \cite{barrett2010inertia}, $\mr(G)=\mr(G-v)$ and it is easy to see that $Z(G) \le Z(G-v)+1.$  Then
\begin{align*}
M(G)&=n-\mr(G)\\ &=1+(n-1)-\mr(G-v)\\ &=1+M(G-v)\\  &=1+Z(G-v)\ge Z(G).
\end{align*}
\end{proof}

\Cref{domvertexprop} suggests a more general technique we term lifting, which we discuss in \Cref{liftingsec}.

Finally, there is a classification for graphs $G$ with $\mr(G)=2$, which is especially useful when determining the minimum rank of small graphs.

\begin{theorem}\cite{barrett2004graphs}\label{minrank29v} Let $G$ be a connected graph on fewer than 9 vertices.  Then $\mr(G) \le 2$ if and only if $G$ does not contain as an induced subgraph any of the graphs in $\Cref{fim}$.
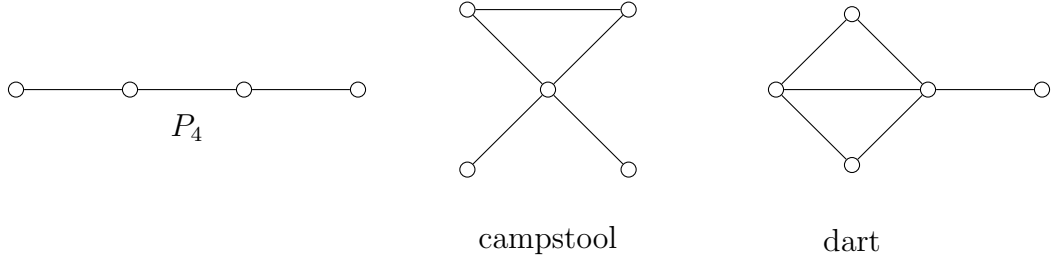
\begin{figure}
    \centering

\begin{tikzpicture}[scale=1]
    \tikzstyle{every node}=[minimum width=0pt, inner sep=2pt, circle]

    \node[draw] (1) at (0,0) {};
    \node[draw] (2) at (1.5,0) {};
    \node[draw] (3) at (3,0) {};
    \node[draw] (4) at (4.5,0) {};

    \draw (1)--(2)--(3)--(4);
    \node at (2.25,-.5) {$P_4$};

    \begin{scope}[shift={(7,0)}]
        % \foreach \n/\i in {1/45, 2/135, 3,225, 4/315}{
        % \node[draw] (\i) at (\i:2cm) {$\n$};
        % }
        \def\n{4}
        % Calculate the angle between vertices
        \foreach \i in {1,...,\n} {    
        % Calculate the coordinates for each vertex    
        \node[circle, draw, minimum size=5pt, inner sep=2pt] (v\i) at ({360/\n*\i +45} : 1.5) {};}
        \node[draw] (5) at (0,0) {};
        \draw (v2)--(5)--(v4)--(v1)--(5) -- (v3);
        \node (a) at (0,-2) {\text{campstool}};
    \end{scope}

    \begin{scope}[shift={(11,0)}]
        \def\n{4}
        % Calculate the angle between vertices
        \foreach \i in {1,...,\n} {    
        % Calculate the coordinates for each vertex    
        \node[circle, draw, minimum size=5pt, inner sep=2pt] (v\i) at ({360/\n*\i } : 1) {};}
        \node[draw] (5) at (2.5,0) {};
        \draw (5) -- (v4) -- (v3)--(v2)--(v1)--(v4)--(v2);
                \node (a) at (0,-2) {\text{dart}};

    \end{scope}
    
\end{tikzpicture}
    \caption{Forbidden induced subgraph obstructions for graphs of minimum rank at most two.}
    \label{fim}
\end{figure}

% \raisebox{.5in}{P4 \includegraphics{P4}}  campstool \includegraphics[width=1in]{foldingstool}
% dart \includegraphics{dart}

\end{theorem}

\section{Partial $k$-Paths}\label{kpathssec}
The class of graphs known as $k$-trees has been studied previously in the context of the minimum rank of a graph.  See for example \cite{vancomppartialktrees,van3connmaxnul3}.  

\begin{definition}\cite{markenzon,bickletreesurvey}
    The complete graph on $k+1$ vertices, $K_{k+1}$, is a $k$-tree.  A $k$-tree with $n$ vertices, $n>k+1$ can be constructed from a $k$-tree with $n-1$ vertices by adding a vertex adjacent to all vertices of a $k$-clique of the existing $k$-tree, and only to these vertices.
\end{definition}

\begin{definition}\cite{markenzon,bickletreesurvey}
    A $k$-path graph $G$ is an alternating sequence of distinct $k$- and $k+1$-cliques $e_0, t_1, e_1, t_2,\ldots, t_p, e_p$, $p\geq 1$, starting and ending with a $k$-clique and such that $t_i$ contains exactly two $k$-cliques $e_{i-1}$ and $e_i$ for $1\leq i \leq p$.
\end{definition}

\begin{lemma}\cite{markenzon,bickletreesurvey}
    A $k$-path is a $k$-tree with exactly two vertices of degree $k$.
\end{lemma}

\begin{example}\label{eg:3path}
    In the graph in \Cref{fig:3path} the sequence of distinct $3$- and $4$-cliques are as follows:
\begin{figure}
    \centering
    \begin{tikzpicture}[scale=2,thick]
		\tikzstyle{every node}=[minimum width=0pt, inner sep=2pt, circle]
			\draw (-3.35,1.27) node[draw] (0) {0};
			\draw (-2.94,1.8) node[draw] (1) {1};
			\draw (-2.72,0.79) node[draw] (2) {2};
			\draw (-2.07,1.23) node[draw] (3) {3};
			\draw (-1.98,2.03) node[draw] (4) {4};
			\draw (-1.57,0.55) node[draw] (5) {5};
			\draw (-1.17,1.47) node[draw] (6) {6};
			\draw (-0.7,0.66) node[draw] (7) {7};
			\draw (-0.41,2.02) node[draw] (8) {8};
			\draw (0.49,1) node[draw] (9) {9};
			\draw (0.38,1.93) node[draw] (10) {10};
			\draw  (0) edge (1);
			\draw  (1) edge (2);
			\draw  (0) edge (2);
			\draw  (2) edge (3);
			\draw  (1) edge (3);
			\draw  (0) edge (3);
			\draw  (1) edge (4);
			\draw  (3) edge (4);
			\draw  (0) edge (4);
			\draw  (4) edge (6);
			\draw  (3) edge (6);
			\draw  (1) edge (6);
			\draw  (4) edge (5);
			\draw  (3) edge (5);
			\draw  (5) edge (6);
			\draw  (3) edge (7);
			\draw  (6) edge (7);
			\draw  (5) edge (7);
			\draw  (6) edge (9);
			\draw  (3) edge (9);
			\draw  (7) edge (9);
			\draw  (6) edge (8);
			\draw  (8) edge (9);
			\draw  (7) edge (8);
			\draw  (8) edge (10);
			\draw  (9) edge (10);
			\draw  (7) edge (10);
    \end{tikzpicture}
    \caption{The 3-path used in \Cref{eg:3path}.}
    \label{fig:3path}
\end{figure}
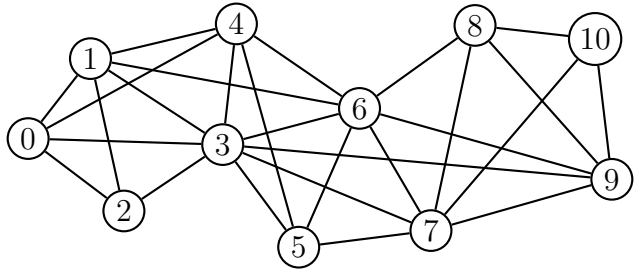

   $ \begin{array}{cc}
    e_0 = \{0,1,2\} & t_1 = \{0,1,2,3\} \\
    e_1 = \{0,1,3\} & t_2 = \{0,1,3,4\} \\
    e_2 = \{1,3,4\} & t_3 = \{1,3,4,6\} \\
    e_3 = \{3,4,6\} & t_4 = \{3,4,5,6\} \\
    e_4 = \{3,5,6\} & t_5 = \{3,5,6,7\} \\
    e_5 = \{3,6,7\} & t_6 = \{3,6,7,9\} \\
    e_6 = \{6,7,9\} & t_7 = \{6,7,8,9\} \\
    e_7 = \{7,8,9\} & t_8 = \{7,8,9,10\} \\
    e_8 = \{8,9,10\} & \\
    \end{array} $

\medskip

Starting with an initial colored set of $e_0$, the sequence of distinct $3$- and $4$-cliques induces a set of zero forcing paths.  In general, $t_i\backslash e_i$ forces $t_i\backslash e_{i-1}$ for $1\leq i \leq 8$.

\medskip

$\begin{array}{cc}
t_1 \backslash e_1 \rightarrow t_1 \backslash e_0 & 2 \rightarrow 3 \\ 
t_2\backslash e_2 \rightarrow t_2\backslash e_1 & 0 \rightarrow 4 \\ 
t_3\backslash e_3 \rightarrow t_3\backslash e_2 & 1 \rightarrow 6 \\
t_4\backslash e_4 \rightarrow t_4\backslash e_3 & 4 \rightarrow 5 \\
t_5\backslash e_5 \rightarrow t_5\backslash e_4 & 5 \rightarrow 7 \\
t_6\backslash e_6 \rightarrow t_6\backslash e_5 & 3 \rightarrow 9 \\
t_7\backslash e_7 \rightarrow t_7\backslash e_6 & 6 \rightarrow 8 \\
t_8\backslash e_8 \rightarrow t_8\backslash e_7 & 7 \rightarrow 10 \\
\end{array}$

\medskip

Thus, the three zero forcing chains are $$0\rightarrow 4\rightarrow 5\rightarrow 7\rightarrow 10, 2\rightarrow 3\rightarrow 9 \text{ and } 1\rightarrow 6\rightarrow 8.$$
\end{example}

\begin{lemma}\cite{bickletreesurvey}
    Let $G$ be a $k$-tree.  Then $G$ is $k$-connected.
\end{lemma}

\begin{theorem}\label{kpathMZ}
    Let $G$ be a $k$-path with sequence $e_0, t_1, e_1, t_2, \ldots, t_p, e_p$ of distinct $k$- and $k+1$-cliques.  Then $e_0$ is a zero forcing set for $G$ and the paths generated by the recursive relation $t_i\backslash e_i$ forces $t_i\backslash e_{i-1}$ for $1 \leq i \leq p$ are a set of zero forcing chains for $G$.  Consequently, $Z(G)=M(G)=k$.
\end{theorem}

\begin{proof} Proceed by induction on the order of $G$. The smallest $k$-path is $K_{k+1}$ and $Z(K_{k+1})=M(K_{k+1})=k$.  Assume that the statement of the theorem is true for all $k$-paths of order at least $n\geq k+1$.  Let $G$ be a $k$-path of order $n+1$ with sequence $e_0, t_1, e_1, t_2, \ldots, t_p, e_p$ of distinct $k$- and $k+1$-cliques.  Let $v=t_p\backslash e_{p-1}$.  The vertex $v$  is one of two vertices of degree $k$ as $v$ is adjacent to only the $k$ vertices in $e_{p-1}$.  The graph $G-v$ is a $k$-path on $n$ vertices with sequence of distinct $k$- and $k+1$-cliques $e_0, t_1, e_1, t_2, \ldots t_{p-1}, e_{p-1}$.  By the inductive hypothesis, $e_0$ is a zero forcing set for $G-v$.  The paths created by $t_i\backslash e_{i}$ forces $t_i\backslash e_{i-1}$,  for $1\leq i\leq p-1$, are zero forcing chains for $G-v$.  None of the vertices of $e_{p-1}$ are used to force other vertices of $G-v$. Since these are the only vertices adjacent to $v$ in $G$, the zero forcing chains for $G-v$ are also valid in $G$. Thus, it only remains to force the vertex $v$ in $G$.  The vertex of $e_{p-1}$ that forces $v$ is $t_p\backslash e_p$.  Thus $e_0$ is a zero forcing set for $G$ and $Z(G) \le |e_o| = k.$ Since $G$ is $k$-connected, $Z(G) \ge M(G) \ge \kappa(G)=k$ by Theorem \ref{kappa_xi_M_Z}.  Therefore $Z(G)=M(G)=k$.
    \end{proof}

\begin{lemma}
    Let $G$ be a $k$-path with sequence $e_0, t_1, e_1, t_2, \ldots, t_p, e_p$ of $k$- and $k+1$-cliques.  Let $P$ be the set of zero forcing chains induced by the sequence of $k$- and $k+1$-cliques.  Let $e$ be an edge of one of the chains in $P$.  Then $G-e$ is not $k$-connected.
\end{lemma}

\begin{proof}
    Let $G$ be as described in the theorem. Let $e$ be an edge of one of the zero forcing chains.  Since the set of zero forcing chains was induced by the sequence of $k$- and $k+1$-cliques, there exists an $r$, $1\leq r\leq p$, such that $e=(t_r\backslash e_r, t_r\backslash e_{r-1})$. Let $u=t_r\backslash e_r$ and $w=t_r\backslash e_{r-1}$. 
    By construction, there exists a set $S$ of $k-1$ vertices such that $e_r=w\bigcup S$ and $e_{r-1}=u\bigcup S$. Since $u$ is in $e_{r-1}$ but not in $e_{r}$, $u$ is not in $e_j$ for $r\leq j\leq p$.  Thus, $u$ is only adjacent to vertices in $e_j$ for $1\leq j\leq r$.  Similarly, $w$ is only adjacent to vertices in $e_j$ for $r\leq j \leq p$.  Consider the graph $H$ obtained from $G-e$ by deleting the vertices in $S$.  Since all the vertices in $e_r$ have been deleted with the exception of $w$, and the edge $uw$ has been deleted, $u$ and $w$ are in distinct components.  Thus, the vertex connectivity of $G-e$ is less than or equal $k-1$.  Therefore $G-e$ is not $k$-connected.
\end{proof}
\begin{definition} 
    If $G$ is a $k$-tree ($k$-path), then any subgraph of $G$ is a partial $k$-tree (partial $k$-path).
\end{definition}

\begin{lemma}\cite{van3connmaxnul3,markenzon}\label{completepartialpath}
    If $G$ is a partial $k$-path, then $G$ can be completed to a $k$-path on the same number of vertices.
\end{lemma}

\begin{theorem}\label{partkpath}
    If $G$ is a $k$-connected partial $k$-path, then $Z(G)=M(G)=k$.
\end{theorem}

\begin{proof}
    Let $G$ be as stated in the theorem.  By Lemma \ref{completepartialpath}, complete $G$ to a $k$-path $H$. By Theorem \ref{kpathMZ}, $Z(H)=k$ and the sequence of distinct $k$- and $k+1$-cliques induce a set of $k$ zero forcing chains.  Since $G$ is $k$-connected, all of the edges used in the zero forcing chains for $H$ are present in $G$. Thus, the same zero forcing chains for $H$ can be used to color the vertices of $G$.  Thus $Z(G)\leq k$.  As $G$ is $k$-connected, $Z(G) \ge M(G) \ge k$ by Theorem \ref{kappa_xi_M_Z}.  Therefore, $Z(G)=M(G)=k$.
\end{proof}

\begin{theorem}\label{3connZ4}
    Let $G$ be a 3-connected graph such that $Z(G)=4$.  Then $M(G)=4$.
\end{theorem}

\begin{proof}

    By \Cref{kappa_xi_M_Z}, $4=Z(G)\ge M(G) \ge \xi(G) \ge \kappa(G)  \ge  3$.  If $M(G)=3$, then by Theorem 14 \cite{van3connmaxnul3}, (ii) $\Rightarrow$ (i), $G$ is a partial 3-path.  Then by \Cref{partkpath}, $Z(G)=3$, contradicting the hypothesis.  Therefore, $M(G)=4$.
\end{proof}

\section{Maximum Nullity and Zero Forcing Number Coincide for seven-vertex graphs.} \label{x7vertices}

 Before proceeding to the proof of our main theorem, we first show how the proof of the known result $M(G)=Z(G)$ for all seven-vertex graphs can be simplified considerably using the results above. It suffices to show that equality holds for the 853 connected seven-vertex graphs. We establish the equality by examining various classes of graphs.  \medskip
\begin{enumerate}
    \item Graphs with a cut-vertex:
    \begin{itemize}
        \item[] Calculating $M(G)$ for all seven-vertex graphs with a cut-vertex by (\ref{cutvertM}) and calculating each $Z(G)$ by an exhaustive search, we find that M(G)=Z(G) for each.
    \end{itemize}
    \item Graphs with a $2-$separation:
    \begin{itemize}
        \item[] As in the previous case, one must calculate $M(G)$ for all seven-vertex connected graphs $G$ with a $2-$separation using \Cref{2sepM2}, and also calculate $Z(G)$ for each of these to determine that each pair of these are equal.
    \end{itemize}
    \item $3-$connected graphs:
    \begin{itemize}
        \item[] It then suffices to determine whether or not $M(G)=Z(G)$ for the 136 $3-$connected graphs on 7 vertices.\medskip

If $G$ has a dominating vertex, then $M(G-v)=Z(G-v)$ by Theorem \ref{nleq6}, so $M(G)=Z(G)$ by Proposition \ref{domvertexprop}. \medskip

If $\kappa(G)=Z(G)$, then $M(G)=Z(G)$ by \Cref{kappa_xi_M_Z}.\medskip

Suppose then that $3 \le \kappa(G) < Z(G) \le 6$.\medskip

If $Z(G)=6$, then $G=K_7$ and $M(K_7)=6=Z(K_7)$.\medskip

So we can assume that $\kappa(G) < Z(G) \le 5.$\medskip

If $Z(G) =5$ then $M(G) \le 5$ and $\mr(G) \ge 2$.  We now apply \Cref{minrank29v}.  Checking all the remaining graphs with $Z(G)=5$, we find that none of these contain $P_4$, campstool, or dart.  Then for each, $2 \le \mr(G) \le 2$ which implies that $M(G)=5=Z(G)$.\medskip

The last case is $Z(G)=4.$  By \Cref{3connZ4}, $M(G)=4$.\medskip

So in all cases, $M(G)=Z(G)$. 
    \end{itemize}
\end{enumerate}

\section{Graphs on eight vertices with a cut vertex or a $2-$separation.}\label{CutVertex2Sep}\medskip

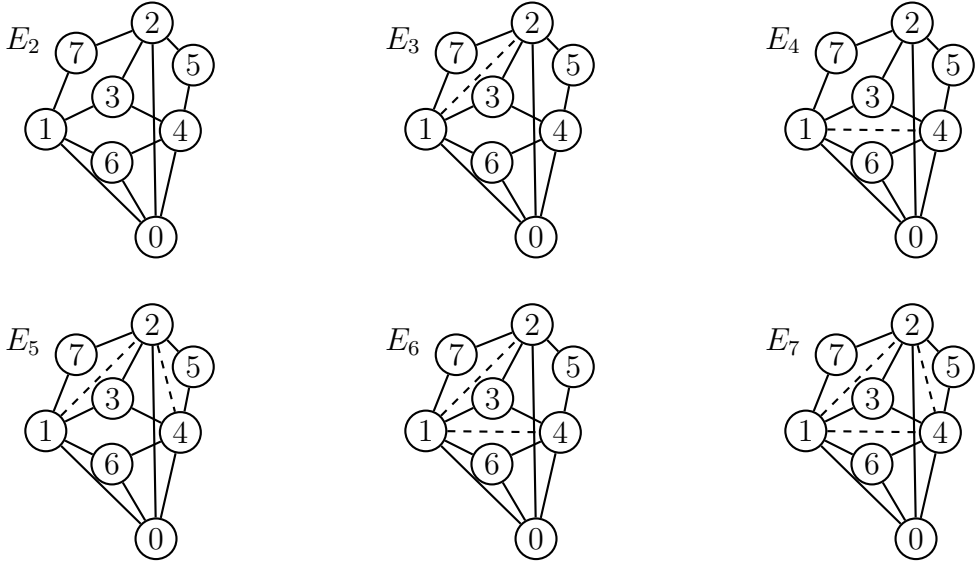
\begin{figure}
    \centering

\begin{tikzpicture}[scale=1.0,thick]
    \tikzstyle{every node}=[minimum width=0pt, inner sep=2pt, circle]

    \node at (-.25,2) (E2) {$E_2$};
	\draw (1.5,-0.61) node[draw] (0) {$0$};
	\draw (0.05,0.82) node[draw] (1) {$1$};
	\draw (1.45,2.23) node[draw] (2) {$2$};
	\draw (0.93,1.25) node[draw] (3) {$3$};
	\draw (1.82,0.8) node[draw] (4) {$4$};
	\draw (1.99,1.67) node[draw] (5) {$5$};
	\draw (0.92,0.38) node[draw] (6) {$6$};
	\draw (0.45,1.83) node[draw] (7) {$7$};
	\draw  (4) edge (6);
	\draw  (0) edge (6);
	\draw  (1) edge (6);
	\draw  (0) edge (2);
	\draw  (2) edge (3);
	\draw  (1) edge (3);
	\draw  (3) edge (4);
	\draw  (0) edge (4);
	\draw  (0) edge (1);
	\draw  (1) edge (7);
	\draw  (2) edge (7);
	\draw  (2) edge (5);
	\draw  (4) edge (5);
% \end{tikzpicture}

% % Graph6: GzcmR?
% $E_3$ \begin{tikzpicture}[scale=1.5,thick]
%     \tikzstyle{every node}=[minimum width=0pt, inner sep=2pt, circle]
\begin{scope}[shift={(5,0)}]
        \node at (-.25,2) (E2) {$E_3$};

	\draw (1.5,-0.61) node[draw] (0) {$0$};
	\draw (0.05,0.82) node[draw] (1) {$1$};
	\draw (1.45,2.23) node[draw] (2) {$2$};
	\draw (0.93,1.25) node[draw] (3) {$3$};
	\draw (1.82,0.8) node[draw] (4) {$4$};
	\draw (1.99,1.67) node[draw] (5) {$5$};
	\draw (0.92,0.38) node[draw] (6) {$6$};
	\draw (0.45,1.83) node[draw] (7) {$7$};
	\draw  (2) edge (5);
	\draw  (4) edge (5);
	\draw  (2) edge (3);
	\draw  (2) edge (7);
	\draw [dashed] (1) edge (2);
	\draw  (0) edge (2);
	\draw  (1) edge (7);
	\draw  (1) edge (3);
	\draw  (3) edge (4);
	\draw  (0) edge (4);
	\draw  (4) edge (6);
	\draw  (0) edge (6);
	\draw  (1) edge (6);
	\draw  (0) edge (1);

 \end{scope}
% \end{tikzpicture}
% %graph6: GrsmR?
% \hskip 1in $E_4$ \begin{tikzpicture}[scale=1.5,thick]
%     \tikzstyle{every node}=[minimum width=0pt, inner sep=2pt, circle]
\begin{scope}[shift={(10,0)}]
        \node at (-.25,2) (E2) {$E_4$};
	\draw (1.5,-0.61) node[draw] (0) {$0$};
	\draw (0.05,0.82) node[draw] (1) {$1$};
	\draw (1.45,2.23) node[draw] (2) {$2$};
	\draw (0.93,1.25) node[draw] (3) {$3$};
	\draw (1.82,0.8) node[draw] (4) {$4$};
	\draw (1.99,1.67) node[draw] (5) {$5$};
	\draw (0.92,0.38) node[draw] (6) {$6$};
	\draw (0.45,1.83) node[draw] (7) {$7$};
	\draw  (0) edge (1);
	\draw  (0) edge (6);
	\draw  (0) edge (4);
	\draw  (0) edge (2);
	\draw  (1) edge (6);
	\draw  (4) edge (6);
	\draw [dashed] (1) edge (4);
	\draw  (1) edge (7);
	\draw  (2) edge (7);
	\draw  (2) edge (5);
	\draw  (4) edge (5);
	\draw  (2) edge (3);
	\draw  (3) edge (4);
	\draw  (1) edge (3);
 \end{scope}
% \end{tikzpicture}

% \vskip 1cm

% %graph6: GzkmR?
% \qquad $E_5$\begin{tikzpicture}[scale=1.5,thick]
%     \tikzstyle{every node}=[minimum width=0pt, inner sep=2pt, circle]

\begin{scope}[shift={(0,-4)}]
        \node at (-.25,2) (E2) {$E_5$};
	\draw (1.5,-0.61) node[draw] (0) {$0$};
	\draw (0.05,0.82) node[draw] (1) {$1$};
	\draw (1.45,2.23) node[draw] (2) {$2$};
	\draw (0.93,1.25) node[draw] (3) {$3$};
	\draw (1.82,0.8) node[draw] (4) {$4$};
	\draw (1.99,1.67) node[draw] (5) {$5$};
	\draw (0.92,0.38) node[draw] (6) {$6$};
	\draw (0.45,1.83) node[draw] (7) {$7$};
	\draw  (4) edge (5);
	\draw  (2) edge (5);
	\draw  (3) edge (4);
	\draw  (0) edge (4);
	\draw  (0) edge (6);
	\draw  (4) edge (6);
	\draw  (0) edge (2);
	\draw [dashed] (1) edge (2);
	\draw  (1) edge (3);
	\draw  (2) edge (3);
	\draw  (2) edge (7);
	\draw  (1) edge (7);
	\draw  (0) edge (1);
	\draw  (1) edge (6);
	\draw [dashed] (2) edge (4);
  \end{scope}
% \end{tikzpicture}
% %graph6: GzsmR?
% \hskip 1in $E_6$
% \begin{tikzpicture}[scale=1.5,thick]
%    \tikzstyle{every node}=[minimum width=0pt, inner sep=2pt, circle]

\begin{scope}[shift={(5,-4)}]
        \node at (-.25,2) (E2) {$E_6$};
	\draw (1.5,-0.61) node[draw] (0) {$0$};
	\draw (0.05,0.82) node[draw] (1) {$1$};
	\draw (1.45,2.23) node[draw] (2) {$2$};
	\draw (0.93,1.25) node[draw] (3) {$3$};
	\draw (1.82,0.8) node[draw] (4) {$4$};
	\draw (1.99,1.67) node[draw] (5) {$5$};
	\draw (0.92,0.38) node[draw] (6) {$6$};
	\draw (0.45,1.83) node[draw] (7) {$7$};
	\draw  (0) edge (6);
	\draw  (0) edge (4);
	\draw  (0) edge (1);
	\draw  (1) edge (6);
	\draw  (0) edge (2);
	\draw [dashed] (1) edge (2);
	\draw  (4) edge (6);
	\draw  (4) edge (5);
	\draw  (2) edge (5);
	\draw  (2) edge (3);
	\draw  (3) edge (4);
	\draw  (1) edge (3);
	\draw  (1) edge (7);
	\draw  (2) edge (7);
	\draw [dashed] (1) edge (4);
  \end{scope}
% \end{tikzpicture}

% \vskip 1cm
% %graph6: Gz{mR
% \quad $E7$
% \begin{tikzpicture}[scale=1.5,thick]
%     \tikzstyle{every node}=[minimum width=0pt, inner sep=2pt, circle]
\begin{scope}[shift={(10,-4)}]
        \node at (-.25,2) (E2) {$E_7$};
	\draw (1.5,-0.61) node[draw] (0) {$0$};
	\draw (0.05,0.82) node[draw] (1) {$1$};
	\draw (1.45,2.23) node[draw] (2) {$2$};
	\draw (0.93,1.25) node[draw] (3) {$3$};
	\draw (1.82,0.8) node[draw] (4) {$4$};
	\draw (1.99,1.67) node[draw] (5) {$5$};
	\draw (0.92,0.38) node[draw] (6) {$6$};
	\draw (0.45,1.83) node[draw] (7) {$7$};
	\draw  (0) edge (6);
	\draw  (0) edge (1);
	\draw [dashed] (1) edge (4);
	\draw  (0) edge (4);
	\draw  (4) edge (6);
	\draw  (4) edge (5);
	\draw  (2) edge (5);
	\draw [dashed] (2) edge (4);
	\draw [dashed] (1) edge (2);
	\draw  (3) edge (4);
	\draw  (2) edge (3);
	\draw  (1) edge (3);
	\draw  (2) edge (7);
	\draw  (1) edge (7);
	\draw  (0) edge (2);
	\draw  (1) edge (6);
 \end{scope}
\end{tikzpicture}
    \caption{The graphs $E_2, E_3, \dots , E_7$}
    \label{fig:E2E7}
\end{figure}

In light of \Cref{cutvertM} and \Cref{2sepM2}, it is natural to begin a search for graphs $G$ with the property $M(G)<Z(G)$ among graphs that possess a cut-vertex or 2-separation.
\begin{enumerate}
    \item  Eight-vertex graphs with a cut-vertex.\medskip

We run through the graphs on 8 vertices with a cut-vertex calculating $M(G)$ by \Cref{cutvertM} and comparing each with $Z(G)$. We find that $M(G)=Z(G)$ for all but one of these graphs, the graph $E_1$ in \Cref{fig:E1H7}.  For $E_1$ we have $M(E_1)=2$ and $Z(E_1)=3$.  This graph is well-known in multiple contexts; for example, it is the simplest of the exceptional graphs given in Appendix B of \cite{johnson2009graphs}.  As this is the first exceptional graph so far we give an argument that $Z(E_1)=3$ and $M(E_1)=2$. It is obvious that the set of three pendent vertices is a zero forcing set, so $Z(E_1) \le 3.$  Examination of all two element subsets of $V(E_1)$ shows that none is a zero forcing set, so $Z(E_1)=3$.\medskip

\begin{figure}
    \centering
    \begin{tikzpicture}[scale=2,thick]
\tikzstyle{every node}=[minimum width=12pt, inner sep=2pt, circle]
    \draw (-0.59,2.35) node[draw] (0) {  };
    \draw (-0.59,1.69) node[draw] (1) { };
    \draw (-1.29,1.69) node[draw] (2) {};
    \draw (0.09,1.69) node[draw] (3) {};
    \draw (-0.92,1.06) node[draw] (4) {};
    \draw (-0.27,1.06) node[draw] (5) {};
    \draw (-1.29,0.49) node[draw] (6) {};
    \draw (0.09,0.49) node[draw] (7) {};
    \draw  (0) edge (1);
    \draw  (1) edge (3);
    \draw  (1) edge (2);
    \draw  (1) edge (4);
    \draw  (1) edge (5);
    \draw  (3) edge (5);
    \draw  (4) edge (5);
    \draw  (2) edge (4);
    \draw  (4) edge (6);
    \draw  (5) edge (7);
    \node (a) at (-.59,.2) {$E_1$};

    \begin{scope}[shift={(3,0)}]
        \draw (-0.59,1.69) node[draw] (1) {v};
    \draw (-1.29,1.69) node[draw] (2) {};
    \draw (0.09,1.69) node[draw] (3) {};
    \draw (-0.92,1.06) node[draw] (4) {};
    \draw (-0.27,1.06) node[draw] (5) {};
    \draw (-1.29,0.49) node[draw] (6) {};
    \draw (0.09,0.49) node[draw] (7) {};
    \draw  (1) edge (3);
    \draw  (1) edge (2);
    \draw  (1) edge (4);
    \draw  (1) edge (5);
    \draw  (3) edge (5);
    \draw  (4) edge (5);
    \draw  (2) edge (4);
    \draw  (4) edge (6);
    \draw  (5) edge (7);
    \node (a) at (-.59,.2) {$H_7$};
    \end{scope}
\end{tikzpicture}
    \caption{The graphs $E_1$ and $H_7$.}
    \label{fig:E1H7}
\end{figure}
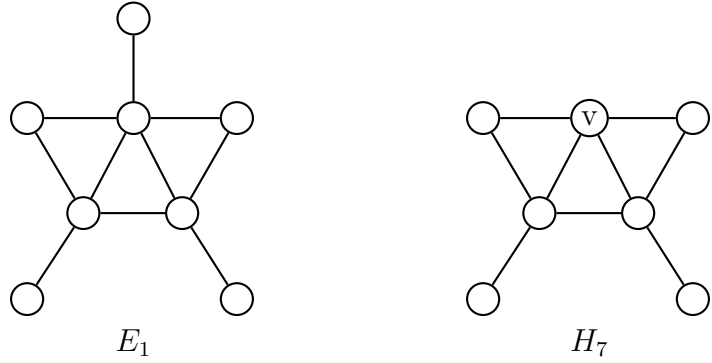

% \begin{figure}
%     \centering
%     \begin{tikzpicture}[scale=2,thick]
% \tikzstyle{every node}=[minimum width=12pt, inner sep=2pt, circle]
    
% \end{tikzpicture}
%     \caption{The graph $H_7$.}
%     \label{fig:H7}
% \end{figure}
% \begin{tikzpicture}[scale=1,thick]
% \tikzstyle{every node}=[minimum width=0pt, inner sep=2pt, circle]
%     \draw (-2.22,1.91) node[draw] (0) {};
%     \draw (-1.82,1.52) node[draw] (1) {};
%     \draw (-2.22,1.15) node[draw] (2) {};
%     \draw (-1.36,1.52) node[draw] (3) {};
%     \draw (-0.96,1.91) node[draw] (4) {};
%     \draw (-0.96,1.15) node[draw] (5) {};
%     \draw  (1) edge (3);
%     \draw  (3) edge (4);
%     \draw  (3) edge (5);
%     \draw  (1) edge (2);
%     \draw  (0) edge (1);
% \end{tikzpicture}

Consider the graph $H_7$ in \Cref{fig:E1H7}, which is graph 482 in \cite{read1998atlas}.  According to the table in \cite{deloss2008table}, $\mr(H_7)=5$. The graph $H_7-v$ is a tree $H$ on 6 vertices and has zero forcing number equal to 2, so $\mr(H) =4$.  Then by the cut-vertex theorem, \[ \mr(E_1)= \min\{\mr(H_7) + \mr(K_2), \mr(H) +\mr(K_2-v)+2\}=\min\{6, 6\} = 6.    \]

Thus $M(E_1)=8-\mr(E_1) = 8 - 6 = 2 < 3 =Z(E_1)$, so $M(E_1) \ne Z(E_1)$.\medskip

\item Eight-vertex graphs with a $2-$separation
\medskip

Running through these graphs and comparing $M(G)$ and $Z(G)$ we find six graphs with $M(G)=3<4=Z(G)$; see \Cref{fig:E2E7}.  The dashed lines denote edges not present in the first graph, $E_2.$
\end{enumerate}

Once the 6 graphs in \Cref{fig:E2E7} have been identified, there is a simpler way to see that $M(E_i)=3$ and $Z(E_i)=4$ for $i=1\dots 7$.  It depends on one more graph theory parameter, defined by a variation of the standard color change rule. We include this because it gives a much clearer picture of why $M(G)$ and $Z(G)$ differ for these graphs. 

Given a simple undirected graph $G$ on $n$ vertices, we consider a graph $\widehat{G}$ with the same vertices and edges, but with loops on some of the vertices.  There are $2^n$ such graphs for each $G$.  Then 

\[S(\widehat{G}) = \{ A \in S(G) \,|\, a_{ii} \ne 0 \ \text{if and only if there is a loop at}\ i \}.\]

The color change rule for graphs with loops is a relaxation of the standard zero forcing color change rule:  

\begin{itemize}

\item[] \textbf{Looped graph color change rule:} If {\em any} vertex has exactly one white neighbor (possibly itself), then that neighbor can be colored blue.  
  
\end{itemize}

The zero forcing number of a loop graph is the minimum number of vertices needed to force every vertex to be blue applying this amended rule.  $\widehat{Z}(G)$ is the maximum of $Z(\widehat{G})$ over all $2^n$ loop graphs.  By \cite{barioli2013parameters}, $M(G) \le \widehat{Z}(G)$.  

We now use $\widehat{Z}(G)$ to simultaneously show that $M(E_i)\le 3, i=2\dots 7$. We only need to consider two cases. 
\begin{enumerate}
\item[a.] There is no loop on vertex 7
\item[b.] There is a loop on vertex 7
\end{enumerate}

Regardless of which of these $E_i$ we choose, the same forcing sets and forcing sequences work for each.

\begin{enumerate}

\item[a.]  Color vertices 0, 1, 3 blue.  Then the sequence of forces 

\[ 7 \rightarrow 2,\ 3 \rightarrow 4,\ 0 \rightarrow 6,\ 4 \rightarrow 5,\ 2\rightarrow 7 \ \]

colors all vertices blue.

\item[b.] Color vertices 3, 4, 5 blue.  Then the sequence 

\[ 5 \rightarrow 2,\ 3 \rightarrow 1,\ 7 \rightarrow 7,\ 2 \rightarrow 0,\ 0 \rightarrow 6 \]

colors all vertices blue. 
     
\end{enumerate}
It follows that $M(E_j) \le \widehat{Z}(E_j) = 3$.  A zero forcing program shows that $Z(G_i)=4, i=2 \dots 7$\medskip

We now recall our main theorem. 

\mainthm*

We have already verified the reverse implication; we now approach the forward implication.

% \section{$3-$connected graphs on 8 vertices}

% \begin{theorem}
% If $G$ is a $3-$connected graph on 8 vertices, then M(G)=Z(G).
% \begin{proof}
% If $G$ has a dominating vertex, then $M(G)=Z(G)$ by Proposition 2.4.\medskip

% It remains to consider graphs with 
%  $3 \le \kappa(G) \le M(G) \le 7$.  \medskip

% If $\kappa(G)=Z(G)$, then M(G)=Z(G), so we can assume that $\kappa(G) < Z(G).$\medskip

% If $Z(G)=7$, then $G=K_8$ and $M(K_8)=7=Z(K_8)$.\medskip

% If $Z(G) =6$ then $M(G) \le 6$ and $\mr(G) \ge 2$.  We now apply Theorem 2.5  Checking all the $3-$connected graphs with $Z(G)=6$, we find that none of these contain $P_4$, campstool, or dart.  Then for each, $2 \le \mr(G) \le 2$ which implies that $M(G)=6=Z(G)$.\medskip

% If $Z(G) = 5$, then $\kappa(G)$ equals 3 or 4. ????????????????\medskip

% If $Z(G)=4$, then by Theorem 3.4, $M(G)=4.$

% \end{proof}
% \end{theorem}

We can construct minimum rank matrices for each graph $E_2,\ldots, E_7$ using the same clique cover.  Notice that $E_7$ can be covered by the following five cliques $V_1=\{1,2,3,4\}, V_2=\{1,2,7\}, V_3= \{2,4,5\}, V_4=\{0,1,4,6\},$ and $V_5=\{0,2\}$.  We can construct a matrix with rank at most $5$ of the form
\renewcommand{\arraystretch}{1.3}
\[ \left[\begin{array}{cccccccc}
a_{4} + a_{5} & a_{4} & a_{5} & 0 & a_{4} & 0 & a_{4} & 0 \\
a_{4} & a_{1} + a_{2} + a_{4} & a_{1} + a_{2} & a_{1} & a_{1} + a_{4} & 0 & a_{4} & a_{2} \\
a_{5} & a_{1} + a_{2} & a_{1} + a_{2} + a_{3} + a_{5} & a_{1} & a_{1} + a_{3} & a_{3} & 0 & a_{2} \\
0 & a_{1} & a_{1} & a_{1} & a_{1} & 0 & 0 & 0 \\
a_{4} & a_{1} + a_{4} & a_{1} + a_{3} & a_{1} & a_{1} + a_{3} + a_{4} & a_{3} & a_{4} & 0 \\
0 & 0 & a_{3} & 0 & a_{3} & a_{3} & 0 & 0 \\
a_{4} & a_{4} & 0 & 0 & a_{4} & 0 & a_{4} & 0 \\
0 & a_{2} & a_{2} & 0 & 0 & 0 & 0 & a_{2}
\end{array}\right] \]
\renewcommand{\arraystretch}{1}
and choose different nonzero values for $a_1,a_2,a_3,a_4,a_5$ to achieve a rank $5$ matrix for each $E_i$, $i \in \{ 2,3, \dots 7\}$, with the associated zero/nonzero pattern.

\section{Lifting minimum rank matrices}\label{liftingsec} 

Since the minimum rank of all graphs on seven vertices is known, \Cref{domvertexprop} suggests a technique to determine witnesses for the minimum rank of graphs $G$ on eight vertices if $G$ has the same minimum rank as one its induced seven-vertex graphs. %{\color{blue}We define the related parameter $\mathfrak{r}_v(G) = n - Z(G) - \mr(G-v)$. It is immediate that $\mathfrak{r}_v(G) \leq r_v(G)$, with equality if and only if $M(G) = Z(G)$.  Note that $M(G)$ can be computed for all such graphs, and thus a matrix $A \in \mathcal{S}(G)$ witnessing $\nullity A = Z(G)$ is a certificate for $\mr(G)$. We describe techniques to determine a witness $A$ with this property in this section.

 \begin{observation}\label{liftlemma}
        Let $G$ be a graph, $v \in V(G)$, and $A \in S(G-v) \subseteq \mathbb{R}^{n-1 \times n-1}$ be a witness for $\mr(G-v)$.  Given $w \in \mathbb{R}^{n-1}$ construct the matrix $B = \begin{bmatrix} A & Aw \\ w^TA & w^TAw \end{bmatrix}$.  Then $\rank B = \rank A$.  Furthermore, if $B \in S(G)$, then $B$ is a witness for $\mr(G)$.
 \end{observation}
 \begin{proof} $\mr(G) \le \rank B = \rank A=\mr(G-v) \le \mr(G).$
\end{proof}

 We term the matrix $B$ constructed in \Cref{liftlemma} a lifting of $A$.  To ensure a lifting is in $S(G)$, we make an observation in \Cref{nullsupp}.  Let $\N(C)$ denote the (right) nullspace of a matrix $C$.  A submatrix $A[\alpha,\beta]$ of a matrix $A$ is the matrix lying in the rows $\alpha$ and the columns $\beta$.  Let $L = \{ 1,2, \dots, n\}$ and $L_v = L \setminus \{v\}$.
 
 %and given a set $S \subseteq V(G)$, we denote by $A[S]$ the submatrix of $A$ consisting of rows corresponding to $S$.

 \begin{lemma}\label{nullsupp}
     Let $B$ be a lifting of $A$, a witness for $\mr(G-v)$.  Then $B$ is a witness of $\mr(G) $ if and only if:
     \begin{enumerate}
         \item $w \in \N(A[L_v \setminus N_G(v), L_v])$
         \item $[Aw]_i \neq 0$ for all $i \in N_G(v)$.
      \end{enumerate}
 \end{lemma}
 \begin{proof}
     The first condition ensures that $B$ has zeroes in all entries corresponding to non-neighbors of $v$ in $G$, while the second ensures that $B$ has nonzero entries corresponding to neighbors of $v$.  %The structure of $B$ guarantees that $B$ is symmetric and the diagonal entry does not contribute to $\rank(B)$.
 \end{proof}

 It is easy to ensure that the first condition in \Cref{nullsupp} is satisfied, but it is not always possible to satisfy the second condition.  For instance, the second condition will never be satisfied if a row of $A[N_G(v),L_v]$ is equal to a row of $A[L_v\setminus N_G(v), L_v].$  For implementation, we adopt a naive approach for choosing $w\in \N(A[L_v\setminus N_G(v), L_v]).$  Our code computes a basis for  $\N(A[L_v\setminus N_G(v), L_v]),$  say $\{ x_1 , x_2, \dots , x_r \}$ and uses $w = \sum_{i=1}^r x_i$. It then determines if $B \in S(G)$, in which case we conclude $r_v(G) =0$ and $Z(G) = M(G)$.

 In practice the naive guess is quite good, but in some cases we needed to be careful with the null vector chosen.  Using this technique, we were able to show that 181 of the remaining 403 graphs at this stage satisfy $M(G) = Z(G)$. 

 \begin{example} \label{lifting_example}
    We illustrate the technique with examples. 
 Consider the graph $G = K_{2,3}$ as labeled in \Cref{house}.  Consider the subgraph $H$ of $G$ obtained by deleting vertex $5$.  We adopt the notation of \Cref{liftlemma}.  Then we have a known minimum rank witnessing matrix for $H$, namely \[A= \begin{bmatrix}
    0&0&1&1 \\ 0&0&1&1 \\ 1&1&0&0 \\ 1&1&0&0
    \end{bmatrix}.\]
    We have $L_v=\{1,2,3,4\}$, $N_G(5) = \{ 1,2\}$ and $L_v\setminus N_G(5)=\{3,4\}$, and thus 
    \[\mathcal{N}(L_5\setminus N_G(5), L_5)=\left\langle(1,-1,0,0)^T, (0,0,1,0)^T, (0,0,0,1)^T\right\rangle.\]
Consider $w = (1,-1,1,1)^T \in \mathcal{N}(A[\{3,4\},\{1,2,3,4\}])$. Then by computing $Aw  = (2,2,0,0)^T$, we see that $w$ satisfies the conditions in \Cref{nullsupp}.  We compute $w^T Aw = x = 0$, and our minimum rank matrix for $G$ becomes
    \[M=\begin{bmatrix}
        0&0&1&1&2 \\ 0&0&1&1&2 \\ 1&1&0&0&0 \\ 1&1&0&0&0 \\ 2&2&0&0&0
    \end{bmatrix}\]

\begin{figure}
    \centering
    \begin{tikzpicture}
          % Define vertices
  \node[circle, draw] (1) at (0,0) {$1$};
  \node[circle, draw] (4) at (2,0) {$4$};
  \node[circle, draw] (3) at (0,2) {$3$};
  \node[circle, draw] (2) at (2,2) {$2$};
    \node[circle, draw] (5) at (1,1) {$5$};

  % Define edges
  \draw (1) -- (4);
  \draw (2) -- (3);
  \draw (2) -- (4);
  \draw (3) -- (1);
  \draw (2) -- (5) -- (1);
    \end{tikzpicture}
    \caption{The graph $G$ in \Cref{lifting_example}.}
    \label{house}
\end{figure}
\end{example}

\begin{example}\label{lifting_nonexample}
    \begin{figure}
    \centering
    \begin{tikzpicture}
          % Define vertices
  \node[circle, draw] (1) at (0,0) {$1$};
  \node[circle, draw] (2) at (2,0) {$2$};
  \node[circle, draw] (3) at (0,2) {$3$};
  \node[circle, draw] (4) at (2,2) {$4$};
    \node[circle, draw] (5) at (4,1) {$5$};

  % Define edges
  \draw (1) -- (2);
  \draw (2) -- (3);
  \draw (3) -- (4);
  \draw (3) -- (1);
  \draw (2) -- (5) -- (4);
    \end{tikzpicture}
    \caption{The graph $G$ in \Cref{lifting_nonexample}.}
    \label{house_adj}
\end{figure}

For the graph in \Cref{house_adj}, a matrix witnessing the minimum rank of $G-\{5\}$ is given by \[A = \begin{bmatrix}
    1&1&1&0\\1&1&1&0\\1&1&2&1\\0&0&1&1
\end{bmatrix}.\]
Since $1 \in L_5 \setminus N_G(5)$, $2 \in N_G(5)$, and $A[\{1\},L_5] = A[\{2\},L_5]$ no linear combination of the columns of $A$ can have a zero entry in the first component but a nonzero entry in the second.

%Since $1 \in V(G) \setminus N_G[\{5\}]$, $2 \in N_G(\{5\})$, and $A[\{1\}] = A[\{2\}]$ no linear combination of the columns of $A$ can have a zero entry in the first component but a nonzero entry in the second.
\end{example}

In practice, the new row and column can be inserted anywhere in $A$ to obtain its lift $B$, which allows us to avoid relabeling the graph.  

\begin{example}\label{lesstrivialeg}
 Consider the graph $G$ in \Cref{lesstrivialfig}.  We delete the red-colored vertex to obtain $G-v$, for which a witness of $\mr(G-v)$ is known, namely:
 \[A= \begin{bmatrix}1 & -1 & 0 & 0& 0 &0  & 2 \\ -1 &1 &0 & 2 &0 & 2&0 \\ 0 &0&0 &2 & 0 &2 &  2 \\ 0& 2 & 2& 2&2&2  & 0 \\ 0 & 0 & 0& 2 & 0& 2 &  2 \\ 0 & 2 & 2 & 2& 2&2 &0  \\ 2&0&2&0&2&0&2  \end{bmatrix}\]
 Applying the lifting technique, we obtain a witness for $\mr(G)$, with the added row and column indicated in red:
 \[ \begin{bmatrix}1 & -1 & 0 & 0& 0 &0 & \textcolor{red!50}{-6} & 2 \\ -1 &1 &0 & 2 &0 & 2&\textcolor{red!50}4&0 \\ 0 &0&0 &2 & 0 &2 & \textcolor{red!50}{-2} & 2 \\ 0& 2 & 2& 2&2&2 &\textcolor{red!50}{10} & 0 \\ 0 & 0 & 0& 2 & 0& 2 & \textcolor{red!50}{-2} & 2 \\ 0 & 2 & 2 & 2& 2&2&\textcolor{red!50}{10} &0 \\\textcolor{red!50}{-6} & \textcolor{red!50}4 &\textcolor{red!50}{-2} & \textcolor{red!50}{10}& \textcolor{red!50}{-2} & \textcolor{red!50}{10} &\textcolor{red!50}{14} & \textcolor{red!50}0 \\ 2&0&2&0&2&0&\textcolor{red!50}0&2  \end{bmatrix}\]

    \begin{figure}
    \centering
\begin{tikzpicture}[scale=0.5]
%\foreach \n/\a in {1/45,2/135,3/225,4/-45} {
        %\node[draw,circle] (\n) at (\a:2) {};
    %}	
    \node[draw,circle] (1) at (45:2.2) {$4$};
    \node[draw,circle] (2) at (135:2.2) {$2$};
    %\node[draw,circle] (3) at (225:2) {$7$};
    \node[draw,circle] (4) at (-45:2.2) {$6$};

    \node[draw,circle,fill=red!50] (3) at (225:2.2) {$7$};
    \draw[thick] (4) -- (1) --(2) -- (3) -- (4) -- (2);
    \draw[thick] (3) --(1);
    \begin{scope}[shift={(2.818,0)}]

    \node[draw,circle] (5) at (45:2.2) {$3$};
    \node[draw,circle] (6) at (-45:2.2) {$5$};

    \draw[thick] (1) -- (6) -- (4) -- (5) -- ( 1);
    \draw[thick] (6) --(1); 
    \end{scope}
    
    \node[draw,circle] (7) at (2.5,4.3) {$8$};
        \node[draw,circle] (8) at (0.1,4.3) {$1$};
        \draw[thick] (2) -- (8);
        \draw[thick,in=30,out=90] (8) -- (3);
        \draw[thick] (6) -- (7) -- (5)--(3);
        \draw[thick,out=60,in=30] (7) -- (8);
        \draw[in=30,out=30,thick] (6) edge[bend left] (3);
\end{tikzpicture}
\caption{The graph in \Cref{lesstrivialeg}.}

\label{lesstrivialfig}
\end{figure}
\end{example}

\section{Algorithmic details}\label{algo-details}

We employ code implemented in Sage to assist with our calculations.  All code is available at \href{https://github.com/hunnellm/maximum-nullity}{https://github.com/hunnellm/maximum-nullity}.

\subsection{Previously implemented bounds}

Our goal is to determine the minimum rank of all graphs with eight vertices, and it is enough to consider the 11,117 connected eight-vertex graphs.  Using the program in \cite{deloss2020techniques}, we are left with 1803 connected graphs whose minimum rank is undetermined by previously implemented techniques.  The program does, however, compute lower and upper bounds for the minimum rank of all graphs. We note that the program establishes again  that $M(E_1) = 2 < 3 = Z(E_1)$.

For a graph $G$, there is previously implemented code to compute the enhanced zero forcing number $\widehat{Z}(G)$.  This code identifies $\widehat{Z}(E_i) = 3 < 4 = Z(E_i)$ for $i \in \{2,3,\dots,7\}$ as before.   This establishes that $M(E_i)<Z(E_i)$ for $i \in \{1,2,\dots, 7\}$.

Our remaining work is to show that $M(G)=Z(G)$ for the remaining 1797 connected graphs with eight vertices.

\subsection{Dominating vertices}
Determining the degree sequence of a graph is computationally efficient, which suggests that any simplification of the minimum rank problem from the degree sequence should be implemented early in the algorithm.

Recall that all graphs $G$ with order at most seven satisfy $M(G){=Z(G)}$.  \Cref{domvertexprop} leads to an observation for graphs on eight vertices.  

\begin{corollary}
Let $G$ be graph on eight vertices.  If $G$ has a dominating vertex, then $M(G) = Z(G)$.    
\end{corollary}

Implementing this corollary establishes that $M(G)=Z(G)$ for an additional 254 remaining graphs.  1543 graphs remain at this stage.

\subsection{Graphs with a 2-separation} 

A formula for the maximum nullity of graphs with $\kappa(G)=2$ was given in \Cref{2sepM2}. Our program implements this formula, and establishes $M(G)=Z(G)$ for a further 948 graphs, which leaves 595 graphs on eight vertices whose minimum rank is undetermined.

\subsection{Vertex connectivity}

From \Cref{kappa_xi_M_Z} it follows that for any graph $G$, $\kappa(G) = Z(G)$ implies $M(G)=Z(G)$.  Of our remaining graphs, 14 satisfy $\kappa(G) = Z(G)$ and thus there 581 graphs on eight vertices whose minimum rank is undetermined.

\subsection{Twin vertex reduction}
In a graph $G$, if $v,w \in V(G)$ have the same open neighborhoods they are twins.  Let $v,w$ be twins in $G$. If $vw \in E(G)$ then $v,w$ are adjacent twins, and are independent twins otherwise.

Suppose $G$ is a graph and $A \in S(G)$ is a minimum rank witness of $G$, i.e. $ \mr(G)=\rank(A)$. For $v \in V(G)$, denote by $a_{vv}$ the diagonal entry of $A$ corresponding to $v$.
\begin{proposition}[\cite{deloss2020techniques}]\label{twinprop}
    Let $G$ be a graph with twin vertices $v,w$
    \begin{enumerate}
        \item If $v,w$ are independent twins and $G-w$ has a minimum rank witness such that $a_{vv}=0$, then $\mr(G) = \mr(G-w)$.
        \item If $v,w$ are adjacent twins and $G-w$ has a minimum rank witness such that $a_{vv}\neq 0$, then $\mr(G) = \mr(G-w)$.
    \end{enumerate}
\end{proposition}

\Cref{twinprop} suggests that for graphs with twins $v,w$, the determination of minimum rank can be reduced to smaller graphs in cases where $a_{vv}$ can be determined to be zero or nonzero. Sufficient conditions for this determination were given in \cite{barrett2013diagonal}, where they appear as corollaries.

\begin{proposition}
    Let $G$ be a graph, $v$ be a vertex of $G$, and $G^{\ell}$ the loop configuration with a loop at $v$ and no other loops or nonloops are specified.  If there exists a set of fewer than $M(G)$ vertices such that starting with these vertices blue, every vertex is eventually colored blue under the standard color change rule with the additional condition that $v$ can color itself if it has no white neighbors, then $a_{vv}=0$ for witnesses of $\mr(G)$.
\end{proposition}

\begin{proposition}
    Let $G$ be a graph, $v$ be a vertex of $G$, and $G^{\ell}$ the loop configuration with no loop at $v$ and no other loops or nonloops are specified.  If there exists a set of fewer than $M(G)$ vertices such that starting with these vertices blue, every vertex is eventually colored blue under the standard color change rule with the additional condition that $v$ can color a white neighbor $u$ if it has no other white neighbors, then $a_{vv} \neq 0$ for witnesses of $\mr(G)$.
\end{proposition}

Note that these are not mutually exclusive, and that in each proposition the color change rule is modified by a special case of the looped color change rule. In some cases, no information about the diagonal entry is gained.  Our code implements these results, and thereby determines that $M(G)=Z(G)$ for 84 of the remaining graphs. Thus there are 497 left to consider.

\subsection{Partial $3$-path reduction}
We determine that 94 of the remaining graphs satisfy the hypotheses of \Cref{3connZ4}, and thus we are left with 403 graphs to consider.

\subsection{Lifting techniques}
We now implement the techniques discussed in \Cref{liftingsec}.  This allows us to find witnesses establishing $M(G) = Z(G)$ for 181 of the remaining graphs.

We started by compiling a list (not necessarily complete) of known witnesses for the minimum rank of graphs on seven vertices.  We then form the set of graphs, denoted by $\mathcal{K}$, corresponding to these matrices using the zero-nonzero pattern of the matrix.

For each of our remaining graphs $G$, we compute the set of induced subgraphs $\mathcal{H} = \{ G-v \ | \ v\in V(G) \}$.  For each $H \in \mathcal{K} \cap \mathcal{H}$, we use the techniques in \Cref{liftingsec} to determine a candidate matrix $A$.  If $A \in \mathcal{S}(G)$ and $\nullity(A) = Z(G)$, we have determined a witness for $M(G)$.

There are 222 graphs left to consider.

\subsection{Colin de Verdiere parameter}\label{cdv}
The Colin de Verdiere parameter $\mu(G)$ is a lower bound on $M(G)$:
\begin{theorem} [\cite{barioli2013parameters}]  For any graph $G$, $\mu(G) \le M(G).$ \end{theorem}

\begin{definition} A \emph{contraction} of an edge in a simple graph $G$ is obtained by identifying the two adjacent vertices of the edge and suppressing any loops or multiple edges formed in the process.  A {\em minor} of a graph $G$ is any graph formed by a series of edge deletions, deletions of isolated vertices, and edge contractions.\end{definition} 

\begin{definition}  The Petersen family, $\cal P$, consists of the graphs in \Cref{fig:P5}.
%\includegraphics[width=4in,angle=270]{PetersenFamily(3).pdf}
%\begin{center}
%\includegraphics[width=5.5in]{PetersenFamily.png}
%\end{center}

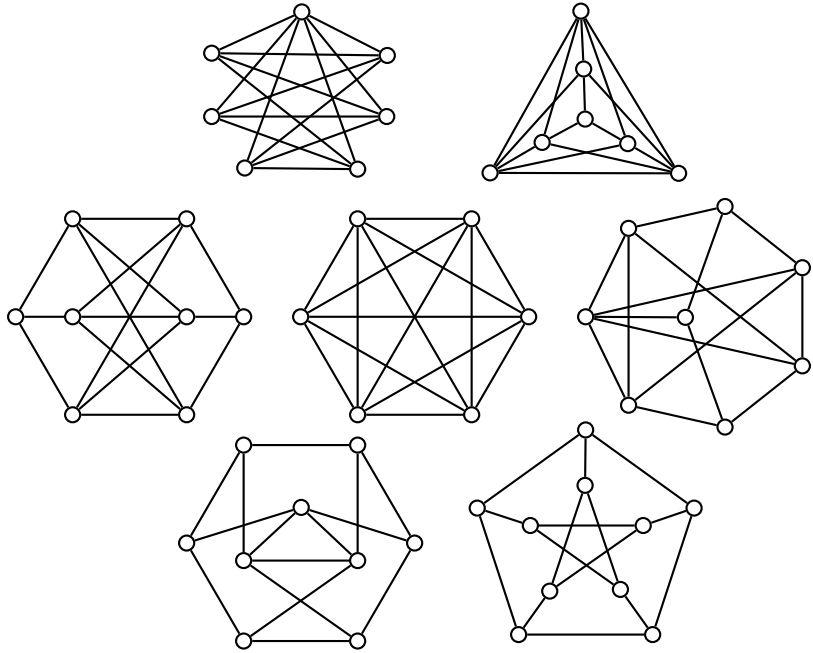
\begin{figure}
    \centering

\begin{tikzpicture}[scale=.75,thick]
    \tikzstyle{every node}=[minimum width=0pt, inner sep=2pt, circle]

\begin{scope}[shift={(5,-4)}]		
        \draw (-2,0) node[draw] (0) {};
        \draw (-1,-1.732050807568877) node[draw] (1) {};
        \draw (1,-1.7320508075688776) node[draw] (2) { };
        \draw (2,0) node[draw] (3) { };
        \draw (1.000000000000001,1.7320508075688767) node[draw] (4) { };
        \draw (-1,1.7320508075688772) node[draw] (5) { };
        \draw  (0) edge (1);
        \draw  (0) edge (2);
        \draw  (0) edge (3);
        \draw  (0) edge (4);
        \draw  (0) edge (5);
        \draw  (1) edge (2);
        \draw  (1) edge (3);
        \draw  (1) edge (4);
        \draw  (1) edge (5);
        \draw  (2) edge (3);
        \draw  (2) edge (4);
        \draw  (2) edge (5);
        \draw  (3) edge (4);
        \draw  (3) edge (5);
        \draw  (4) edge (5);
\end{scope}  
	
\begin{scope}[shift={(8,-8)},rotate=-90]
        \draw (0.82,0.61) node[draw] (0) { };
        \draw (-1.02,-0.01) node[draw] (1) { };
        \draw (0.85,-0.63) node[draw] (2) { };
        \draw (-0.31,1.01) node[draw] (3) { };
        \draw (-0.31,-0.97) node[draw] (4) { };
        \draw (-2,0) node[draw] (5) { };
        \draw (-0.6180339887498951,-1.9021130325903073) node[draw] (6) { };
        \draw (1.618033988749894,-1.1755705045849465) node[draw] (7) { };
        \draw (1.6180339887498951,1.1755705045849458) node[draw] (8) { };
        \draw (-0.6180339887498946,1.9021130325903073) node[draw] (9) { };
        \draw  (0) edge (1);
        \draw  (1) edge (2);
        \draw  (2) edge (3);
        \draw  (3) edge (4);
        \draw  (0) edge (4);
        \draw  (5) edge (6);
        \draw  (6) edge (7);
        \draw  (7) edge (8);
        \draw  (8) edge (9);
        \draw  (5) edge (9);
        \draw  (0) edge (8);
        \draw  (3) edge (9);
        \draw  (1) edge (5);
        \draw  (4) edge (6);
        \draw  (2) edge (7);
\end{scope}

\begin{scope}[shift={(3.2,-.5)}]
        \draw (-0.18,1.89) node[draw] (0) { };
        \draw (-1.76,1.16) node[draw] (1) { };
        \draw (1.32,1.12) node[draw] (2) { };
        \draw (-1.76,0.04) node[draw] (3) { };
        \draw (1.31,0.04) node[draw] (4) { };
        \draw (-1.18,-0.87) node[draw] (5) { };
        \draw (0.8,-0.89) node[draw] (6) { };
        \draw  (0) edge (1);
        \draw  (1) edge (2);
        \draw  (1) edge (4);
        \draw  (1) edge (6);
        \draw  (0) edge (3);
        \draw  (2) edge (3);
        \draw  (3) edge (4);
        \draw  (3) edge (6);
        \draw  (5) edge (6);
        \draw  (4) edge (5);
        \draw  (2) edge (5);
        \draw  (0) edge (5);
        \draw  (0) edge (2);
        \draw  (0) edge (4);
        \draw  (0) edge (6);
\end{scope}

\begin{scope}[shift={(0,-4)}]
        \draw (-2,0) node[draw] (0) { };
        \draw (-1,-1.732050807568877) node[draw] (1) { };
        \draw (1,-1.7320508075688776) node[draw] (2) { };
        \draw (2,0) node[draw] (3) { };
        \draw (1.000000000000001,1.7320508075688767) node[draw] (4) { };
        \draw (-1,1.7320508075688772) node[draw] (5) { };
        \draw (-1,0) node[draw] (6) { };
        \draw (1,0) node[draw] (7) { };
        \draw  (0) edge (1);
        \draw  (1) edge (2);
        \draw  (2) edge (3);
        \draw  (3) edge (4);
        \draw  (4) edge (5);
        \draw  (0) edge (5);
        \draw  (6) edge (7);
        \draw  (2) edge (6);
        \draw  (1) edge (7);
        \draw  (5) edge (7);
        \draw  (4) edge (6);
        \draw  (0) edge (6);
        \draw  (3) edge (7);
        \draw  (2) edge (5);
        \draw  (1) edge (4);
\end{scope}

\begin{scope}[shift={(3,-8)}]
       \draw (-2,0) node[draw] (0) { };
        \draw (-1,-1.732050807568877) node[draw] (1) { };
        \draw (1,-1.7320508075688776) node[draw] (2) { };
        \draw (2,0) node[draw] (3) { };
        \draw (1.000000000000001,1.7320508075688767) node[draw] (4) { };
        \draw (-1,1.7320508075688772) node[draw] (5) { };
        \draw (-1,-0.31) node[draw] (6) { };
        \draw (1,-0.31) node[draw] (7) { };
        \draw (0.01,0.63) node[draw] (8) { };
        \draw  (0) edge (1);
        \draw  (1) edge (2);
        \draw  (2) edge (3);
        \draw  (3) edge (4);
        \draw  (4) edge (5);
        \draw  (0) edge (5);
        \draw  (5) edge (6);
        \draw  (4) edge (7);
        \draw  (6) edge (7);
        \draw  (0) edge (8);
        \draw  (6) edge (8);
        \draw  (7) edge (8);
        \draw  (3) edge (8);
        \draw  (1) edge (7);
        \draw  (2) edge (6);
\end{scope}

\begin{scope}[shift={(10,-4)}]
        \draw (-2,0) node[draw] (0) { };
        \draw (-1.2469796037174672,-1.56366296493606) node[draw] (1) { };
        \draw (0.44504186791262895,-1.9498558243636472) node[draw] (2) { };
        \draw (1.8019377358048383,-0.8677674782351165) node[draw] (3) { };
        \draw (1.8019377358048383,0.8677674782351158) node[draw] (4) { };
        \draw (0.44504186791262895,1.9498558243636472) node[draw] (5) { };
        \draw (-1.2469796037174665,1.5636629649360598) node[draw] (6) { };
        \draw (-0.25,-0.01) node[draw] (7) { };
        \draw  (0) edge (1);
        \draw  (1) edge (2);
        \draw  (2) edge (3);
        \draw  (3) edge (4);
        \draw  (4) edge (5);
        \draw  (5) edge (6);
        \draw  (0) edge (6);
        \draw  (0) edge (7);
        \draw  (0) edge (4);
        \draw  (0) edge (3);
        \draw  (5) edge (7);
        \draw  (2) edge (7);
        \draw  (1) edge (6);
        \draw  (3) edge (6);
        \draw  (1) edge (4);
\end{scope}

\begin{scope}[shift={(8,-.5)},rotate=30]
        \draw (-1.93,0.01) node[draw] (0) { };
        \draw (0.93,-1.65) node[draw] (1) { };
        \draw (0.88,1.69) node[draw] (2) { };
        \draw (-0.87,0.02) node[draw] (3) { };
        \draw (0.42,-0.75) node[draw] (4) {};
        \draw (0.41,0.78) node[draw] (5) { };
        \draw (-0.01,0) node[draw] (6) { };
        \draw  (4) edge (6);
        \draw  (1) edge (4);
        \draw  (5) edge (6);
        \draw  (2) edge (5);
        \draw  (3) edge (6);
        \draw  (0) edge (3);
        \draw  (1) edge (2);
        \draw  (0) edge (1);
        \draw  (0) edge (2);
        \draw  (0) edge (5);
        \draw  (1) edge (5);
        \draw  (2) edge (4);
        \draw  (0) edge (4);
        \draw  (2) edge (3);
        \draw  (1) edge (3);
\end{scope}

\end{tikzpicture}
    \caption{The Petersen family of graphs $\mathcal{P}_5$}
    \label{fig:P5}
\end{figure}

\end{definition}

\begin{theorem}[\cite{schrijver1997minor}, p. 165] The parameter $\mu$ is minor-monotone; i.e., if $H$ is a minor of $G$, then $\mu(H) \le \mu(G).$\end{theorem}

\begin{theorem}[\cite{schrijver1997minor}, p. 188] For a graph $G$, we have $\mu(G) \le 4$ if and only if $G$ does not have a minor in the Petersen family. 
Equivalently, $\mu(G) \ge 5$ if and only if $G$ has a minor in the Petersen family. \end{theorem} 

Of course, in order for $H$ to be a minor of an 8-vertex graph and be a member of the Petersen family, it must be one of the first 5 graphs in the figure above.  Call this group of 5 graphs ${\cal P}_5$.\\  

\begin{theorem}
    If $G$ is an 8-vertex graph with $Z(G)=5$ and $G$ has a minor in ${\cal P}_5$, then $M(G)=Z(G)$.
\end{theorem}

\begin{proof}
  A minor in ${\cal P}_5$ implies $\mu(G)$ is at least five, so $5 \le \mu(G)  \le M(G) \le Z(G) =5$. Thus $M(G)=Z(G)=5$.\\  
\end{proof}  

Implementing this result in our code,  we find 191 additional graphs for which $M(G)=Z(G)$.

\subsection{The final 31 graphs}
\label{final31}
There are $31$ graphs that are not covered by any of the arguments above.  Here we provide a general construction method for minimum rank matrices for these $31$ graphs.  We label each graph with a graph6 string, and provide a minimum rank matrix that corresponds to the vertex labeling of the provided graph6 string.  We would like to thank Tracy Hall for suggesting this construction.

For all graphs $G$ in the last 31 we have $Z(G)=5$, so $\mr(G) \geq 3$.  It is also illuminating to consider the positive semidefinite zero forcing number \cite{barioli2010zero}.  Of the last 31 graphs, 13 of them had $Z_+(G) =4 < 5 = Z(G)$ and we list their graph6 strings here.  

\begin{center}
 $\{$\verb+ Gvd|TO+, \verb+GfF|tW+, \verb+GfD|tW+, \verb+GbD|tW+, \verb+GbD|t[+, \verb+Gvx|Ro+, \verb+Gfx|Ro+, \\
 \verb+GrxX[s+, \verb+GryW[k+, \verb+GrW[[k+, \verb+GvW[[k+, \verb+GrY[[k+, \verb+ Gv][[k+ $\}$
 \end{center}

To verify we have constructed a matrix for each $G$ that has at most rank $3$, we construct $A \in \mathcal{S}(G)$ by choosing a vector representation for the vertices of $G$.  In each case we were able to construct $A=M^TDM$ where $M$ is a $3 \times 8$ matrix containing the column vectors corresponding to the vertices in $G$ and $D= \mathrm{diag}(1,1,-1)$. So in each case we were able to find a minimum rank matrix in $\mathcal{S}(G)$ with two positive eigenvalues and one negative eigenvalue.  

For example we compute $A=M^TDM$ for the graph $G$ with graph6 string \verb+G~xX{s+ which has adjacency matrix

\[ \begin{bmatrix}
                0 & 1 & 1 & 1 & 1 & 0 & 0 & 1 \\
                1 & 0 & 1 & 1 & 1 & 1 & 0 & 0 \\
                1 & 1 & 0 & 1 & 1 & 0 & 1 & 0 \\
                1 & 1 & 1 & 0 & 0 & 1 & 1 & 1 \\
                1 & 1 & 1 & 0 & 0 & 1 & 1 & 1 \\
                0 & 1 & 0 & 1 & 1 & 0 & 1 & 0 \\
                0 & 0 & 1 & 1 & 1 & 1 & 0 & 1 \\
                1 & 0 & 0 & 1 & 1 & 0 & 1 & 0
            \end{bmatrix}. \]

We choose 
\[ M=\begin{bmatrix*}[r]
    1 & 2 & 1 & 1 & 2 & 0 & 2 & 0 \\
0 & 1 & 0 & 1 & 1 & 1 & -4 & 0 \\
2 & 0 & 0 & 1 & 3 & 0 & 1 & 1
\end{bmatrix*} \]
with the columns chosen so that for $i\neq j$, $m_i^T D m_j = 0$ if and only if $ij$ is not an edge for $G$.  

Therefore we have
\[ A = M^T D M= \begin{bmatrix*}[r]
    -3 & 2 & 1 & -1 & -4 & 0 & 0 & -2 \\
2 & 5 & 2 & 3 & 5 & 1 & 0 & 0 \\
1 & 2 & 1 & 1 & 2 & 0 & 2 & 0 \\
-1 & 3 & 1 & 1 & 0 & 1 & -3 & -1 \\
-4 & 5 & 2 & 0 & -4 & 1 & -3 & -3 \\
0 & 1 & 0 & 1 & 1 & 1 & -4 & 0 \\
0 & 0 & 2 & -3 & -3 & -4 & 19 & -1 \\
-2 & 0 & 0 & -1 & -3 & 0 & -1 & -1
\end{bmatrix*} \]
is such that $A\in \mathcal{S}(G)$ and $\rank A=3$. Thus $A$ is a witness for the minimum rank of $G$. The adjacency matrix and a witness for the minimum rank are given for each of the remaining graphs in the appendix.

% \begin{center}
% \begin{tabular}{|c|l|c|c|c|} \hline
%      & Graph6 & $Z_+(G)$ & $\mathrm{Adj}(G)$ & $M^T D M = A$ \\ \hline
%     23 & \verb+'G~xX{s'+ & $5$ & \tiny $\begin{bmatrix*}[r]
%                         0 & 1 & 1 & 1 & 1 & 0 & 0 & 1 \\
%                 1 & 0 & 1 & 1 & 1 & 1 & 0 & 0 \\
%                 1 & 1 & 0 & 1 & 1 & 0 & 1 & 0 \\
%                 1 & 1 & 1 & 0 & 0 & 1 & 1 & 1 \\
%                 1 & 1 & 1 & 0 & 0 & 1 & 1 & 1 \\
%                 0 & 1 & 0 & 1 & 1 & 0 & 1 & 0 \\
%                 0 & 0 & 1 & 1 & 1 & 1 & 0 & 1 \\
%                 1 & 0 & 0 & 1 & 1 & 0 & 1 & 0
%                     \end{bmatrix*}$ & \tiny 
%                     $\begin{bmatrix*}[r]
%                         -3 & 2 & 1 & -1 & -4 & 0 & 0 & -2 \\
%                         2 & 5 & 2 & 3 & 5 & 1 & 0 & 0 \\
%                         1 & 2 & 1 & 1 & 2 & 0 & 2 & 0 \\
%                         -1 & 3 & 1 & 1 & 0 & 1 & -3 & -1 \\
%                         -4 & 5 & 2 & 0 & -4 & 1 & -3 & -3 \\
%                         0 & 1 & 0 & 1 & 1 & 1 & -4 & 0 \\
%                         0 & 0 & 2 & -3 & -3 & -4 & 19 & -1 \\
%                         -2 & 0 & 0 & -1 & -3 & 0 & -1 & -1
%                     \end{bmatrix*}$ \normalsize \TBstrut \\ \hline
                    
% \end{tabular}
% \end{center}

\section{Further observations}

As we investigated commonalities among the graphs $E_i$ for $2\leq i \leq 7$, we noticed that all of these graphs have the same set of zero forcing sets.  This is related to forts \cite{effects2018fast}, another type of distinguished set of vertices related to zero forcing.

\begin{definition}
    For a graph $G$ a subset of vertices $F \subset V$ is a \emph{fort} if $|N_G(v) \cap F| \neq 1$ for all $v \in V \setminus F$.
\end{definition}
Equivalently, a fort is a set of vertices none of which can be zero forced by a vertex outside of the fort.  We denote the set of all forts for a graph $G$ by $\mathrm{Fort}(G)$.
\begin{proposition}[\cite{effects2018fast}]\label{ZF_F_inter}
    For a graph $G$, $B\subset V$ is a zero forcing set if and only if $B\cap F \neq \varnothing$ for all $F \in \mathrm{Fort}(G)$.
\end{proposition}
% \begin{proof}
%     First let $B \subset V$ and $F$ some fort of $G$ such that $B\cap F = \varnothing$.  Then there is no vertex in $V\setminus F \supset B$ that can force any vertex in $F$, so $B$ is not a zero forcing set for $G$.  Therefore if $B$ is a zero forcing set, then $B \cap F \neq \varnothing$ for all $F \in \mathrm{Fort}(G)$.  
    
%     Now assume $B$ is not a zero forcing set.  Let us denote by $\widetilde{B}$ the closure of the set of $B$ with respect to the standard color change rule.  Then $\widetilde{B} \neq V$ so $F = V \setminus \widetilde{B} \in \mathrm{Fort}(G)$ and $B \cap F = \varnothing$.  Therefore if $B \cap F \neq \emptyset$ for all $F \in \mathrm{Fort}(G)$ then $B$ is a zero forcing set for $G$.
% \end{proof}

We denote the set of all zero forcing sets for $G$ by $\mathrm{ZF}(G)$.

\begin{proposition}\label{ZF_Fort}
    For two graphs $G$ and $G'$ we have $\mathrm{ZF}(G) = \mathrm{ZF}(G')$ if and only if $\mathrm{Fort}(G) = \mathrm{Fort}(G')$.
\end{proposition}
\begin{proof}
    Let $G$ and $G'$ be graphs such that $\mathrm{Fort}(G) = \mathrm{Fort}(G')$ and let $B \in \mathrm{ZF}(G)$. Then $B \cap F \neq \varnothing$ for all $F \in \mathrm{Fort}(G)= \mathrm{Fort}(G')$ so by \Cref{ZF_F_inter} we have $B \in \mathrm{ZF}(G')$ so $\mathrm{ZF}(G) \subset \mathrm{ZF}(G')$.  A similar argument shows that $\mathrm{ZF}(G') \subset \mathrm{ZF}(G)$ and we have $\mathrm{ZF}(G) = \mathrm{ZF}(G')$.  

    Now assume that $\mathrm{ZF}(G) = \mathrm{ZF}(G')$.  Let $F \in \mathrm{Fort}(G)$ then $V\setminus F$ contains no zero forcing set from $\mathrm{ZF}(G)$, but $\mathrm{ZF}(G) = \mathrm{ZF}(G')$.  So by \Cref{ZF_F_inter} we have $F \in \mathrm{Fort}(G')$ and similarly for inclusion in the other direction.  Therefore $\mathrm{Fort}(G) = \mathrm{Fort}(G')$.
\end{proof}

In the language of the results above we found that $\mathrm{ZF}(E_i) = \mathrm{ZF}(E_j)$ for all $2 \leq i,j \leq 7$ and so by \Cref{ZF_Fort} we have $\mathrm{Fort}(E_i) = \mathrm{Fort}(E_j)$ for $2 \leq i,j \leq 7$.

\section*{Appendix - Minimum Rank Witnesses}

For the 31 graphs whose minimum rank (maximum nullity) is undetermined by the the preceding methods, the following table compiles witnesses for the minimum rank of the graph.  For a graph $G$, we include two common ways of identifying $G$, namely the $\mathrm{graph6}$-string codes and adjacency matrices, along with a matrix $A\in\mathcal{S}(G)$ which establishes that $M(G)=Z(G)$.
\begin{center}
\begin{longtable}{|c|c|c|c|} \hline
     & Graph6 & Adjacency matrix & $M^T D M = A$ \\ \hline
      1 & \verb+ Gvd|TO + & \small 
$ \left[\begin{array}{rrrrrrrr}
0 & 1 & 1 & 1 & 1 & 0 & 1 & 1 \\
1 & 0 & 0 & 1 & 0 & 1 & 0 & 0 \\
1 & 0 & 0 & 1 & 0 & 1 & 1 & 1 \\
1 & 1 & 1 & 0 & 1 & 1 & 0 & 0 \\
1 & 0 & 0 & 1 & 0 & 1 & 1 & 1 \\
0 & 1 & 1 & 1 & 1 & 0 & 0 & 0 \\
1 & 0 & 1 & 0 & 1 & 0 & 0 & 0 \\
1 & 0 & 1 & 0 & 1 & 0 & 0 & 0
\end{array}\right] $ & \small $ \left[\begin{array}{rrrrrrrr}
-15 & -1 & -4 & -3 & -4 & 0 & -2 & -2 \\
-1 & 1 & 0 & 1 & 0 & 2 & 0 & 0 \\
-4 & 0 & 0 & -1 & 0 & 1 & -1 & -1 \\
-3 & 1 & -1 & 1 & -1 & 2 & 0 & 0 \\
-4 & 0 & 0 & -1 & 0 & 1 & -1 & -1 \\
0 & 2 & 1 & 2 & 1 & 4 & 0 & 0 \\
-2 & 0 & -1 & 0 & -1 & 0 & 0 & 0 \\
-2 & 0 & -1 & 0 & -1 & 0 & 0 & 0
\end{array}\right] $ 
\normalsize \TBstrut \\ \hline

2 & \verb+ GzD~VS + & \small 
$ \left[\begin{array}{rrrrrrrr}
0 & 1 & 1 & 0 & 0 & 0 & 1 & 1 \\
1 & 0 & 1 & 1 & 0 & 1 & 1 & 1 \\
1 & 1 & 0 & 1 & 0 & 1 & 1 & 1 \\
0 & 1 & 1 & 0 & 1 & 1 & 0 & 0 \\
0 & 0 & 0 & 1 & 0 & 1 & 1 & 1 \\
0 & 1 & 1 & 1 & 1 & 0 & 0 & 0 \\
1 & 1 & 1 & 0 & 1 & 0 & 0 & 1 \\
1 & 1 & 1 & 0 & 1 & 0 & 1 & 0
\end{array}\right] $ & \small $ \left[\begin{array}{rrrrrrrr}
29 & 5 & 5 & 0 & 0 & 0 & 11 & 11 \\
5 & 1 & 1 & 1 & 0 & 1 & 2 & 2 \\
5 & 1 & 1 & 1 & 0 & 1 & 2 & 2 \\
0 & 1 & 1 & 5 & 3 & 5 & 0 & 0 \\
0 & 0 & 0 & 3 & -4 & 3 & 1 & 1 \\
0 & 1 & 1 & 5 & 3 & 5 & 0 & 0 \\
11 & 2 & 2 & 0 & 1 & 0 & 4 & 4 \\
11 & 2 & 2 & 0 & 1 & 0 & 4 & 4
\end{array}\right] $ 
\normalsize \TBstrut \\ \hline

3 & \verb+ Gr\~VS + & \small 
$ \left[\begin{array}{rrrrrrrr}
0 & 1 & 1 & 0 & 0 & 0 & 1 & 1 \\
1 & 0 & 0 & 1 & 1 & 1 & 1 & 1 \\
1 & 0 & 0 & 1 & 1 & 1 & 1 & 1 \\
0 & 1 & 1 & 0 & 1 & 1 & 0 & 0 \\
0 & 1 & 1 & 1 & 0 & 1 & 1 & 1 \\
0 & 1 & 1 & 1 & 1 & 0 & 0 & 0 \\
1 & 1 & 1 & 0 & 1 & 0 & 0 & 1 \\
1 & 1 & 1 & 0 & 1 & 0 & 1 & 0
\end{array}\right] $ & \small $ \left[\begin{array}{rrrrrrrr}
1 & 1 & 1 & 0 & 0 & 0 & 2 & 2 \\
1 & 16 & 0 & 2 & 1 & 2 & 9 & 9 \\
1 & 0 & -23 & -9 & -14 & -9 & -1 & -1 \\
0 & 2 & -9 & -3 & -5 & -3 & 0 & 0 \\
0 & 1 & -14 & -5 & -8 & -5 & -1 & -1 \\
0 & 2 & -9 & -3 & -5 & -3 & 0 & 0 \\
2 & 9 & -1 & 0 & -1 & 0 & 7 & 7 \\
2 & 9 & -1 & 0 & -1 & 0 & 7 & 7
\end{array}\right] $ 
\normalsize \TBstrut \\ \hline

4 & \verb+ GfF|tW + & \small 
$ \left[\begin{array}{rrrrrrrr}
0 & 1 & 0 & 1 & 0 & 1 & 1 & 1 \\
1 & 0 & 0 & 1 & 0 & 1 & 0 & 0 \\
0 & 0 & 0 & 1 & 0 & 1 & 1 & 1 \\
1 & 1 & 1 & 0 & 1 & 1 & 1 & 0 \\
0 & 0 & 0 & 1 & 0 & 1 & 1 & 1 \\
1 & 1 & 1 & 1 & 1 & 0 & 0 & 1 \\
1 & 0 & 1 & 1 & 1 & 0 & 0 & 0 \\
1 & 0 & 1 & 0 & 1 & 1 & 0 & 0
\end{array}\right] $ & \small $ \left[\begin{array}{rrrrrrrr}
1 & 1 & 0 & 2 & 0 & 1 & 1 & -1 \\
1 & 1 & 0 & 1 & 0 & 2 & 0 & 0 \\
0 & 0 & 0 & 1 & 0 & -1 & 1 & -1 \\
2 & 1 & 1 & 2 & 1 & 2 & 1 & 0 \\
0 & 0 & 0 & 1 & 0 & -1 & 1 & -1 \\
1 & 2 & -1 & 2 & -1 & 3 & 0 & -1 \\
1 & 0 & 1 & 1 & 1 & 0 & 1 & 0 \\
-1 & 0 & -1 & 0 & -1 & -1 & 0 & -1
\end{array}\right] $ 
\normalsize \TBstrut \\ \hline

5 & \verb+ GfD|tW + & \small 
$ \left[\begin{array}{rrrrrrrr}
0 & 1 & 0 & 1 & 0 & 0 & 1 & 1 \\
1 & 0 & 0 & 1 & 0 & 1 & 0 & 0 \\
0 & 0 & 0 & 1 & 0 & 1 & 1 & 1 \\
1 & 1 & 1 & 0 & 1 & 1 & 1 & 0 \\
0 & 0 & 0 & 1 & 0 & 1 & 1 & 1 \\
0 & 1 & 1 & 1 & 1 & 0 & 0 & 1 \\
1 & 0 & 1 & 1 & 1 & 0 & 0 & 0 \\
1 & 0 & 1 & 0 & 1 & 1 & 0 & 0
\end{array}\right] $ & \small $ \left[\begin{array}{rrrrrrrr}
1 & 1 & 0 & 1 & 0 & 0 & 2 & 2 \\
1 & 1 & 0 & 3 & 0 & 2 & 0 & 0 \\
0 & 0 & 0 & 2 & 0 & 2 & -2 & -2 \\
1 & 3 & 2 & 5 & 2 & 6 & 4 & 0 \\
0 & 0 & 0 & 2 & 0 & 2 & -2 & -2 \\
0 & 2 & 2 & 6 & 2 & 8 & 0 & -4 \\
2 & 0 & -2 & 4 & -2 & 0 & -4 & 0 \\
2 & 0 & -2 & 0 & -2 & -4 & 0 & 4
\end{array}\right] $ 
\normalsize \TBstrut \\ \hline

6 & \verb+ GbD|tW + & \small 
$ \left[\begin{array}{rrrrrrrr}
0 & 1 & 0 & 0 & 0 & 0 & 1 & 1 \\
1 & 0 & 0 & 1 & 0 & 1 & 0 & 0 \\
0 & 0 & 0 & 1 & 0 & 1 & 1 & 1 \\
0 & 1 & 1 & 0 & 1 & 1 & 1 & 0 \\
0 & 0 & 0 & 1 & 0 & 1 & 1 & 1 \\
0 & 1 & 1 & 1 & 1 & 0 & 0 & 1 \\
1 & 0 & 1 & 1 & 1 & 0 & 0 & 0 \\
1 & 0 & 1 & 0 & 1 & 1 & 0 & 0
\end{array}\right] $ & \small $ \left[\begin{array}{rrrrrrrr}
1 & 1 & 0 & 0 & 0 & 0 & 1 & -1 \\
1 & 1 & 0 & 1 & 0 & 1 & 0 & 0 \\
0 & 0 & 0 & -1 & 0 & -1 & 1 & -1 \\
0 & 1 & -1 & 2 & -1 & 1 & -1 & 0 \\
0 & 0 & 0 & -1 & 0 & -1 & 1 & -1 \\
0 & 1 & -1 & 1 & -1 & 0 & 0 & -1 \\
1 & 0 & 1 & -1 & 1 & 0 & 1 & 0 \\
-1 & 0 & -1 & 0 & -1 & -1 & 0 & -1
\end{array}\right] $ 
\normalsize \TBstrut \\ \hline

7 & \verb+ GbD|t[ + & \small 
$ \left[\begin{array}{rrrrrrrr}
0 & 1 & 0 & 0 & 0 & 0 & 1 & 1 \\
1 & 0 & 0 & 1 & 0 & 1 & 0 & 0 \\
0 & 0 & 0 & 1 & 0 & 1 & 1 & 1 \\
0 & 1 & 1 & 0 & 1 & 1 & 1 & 0 \\
0 & 0 & 0 & 1 & 0 & 1 & 1 & 1 \\
0 & 1 & 1 & 1 & 1 & 0 & 0 & 1 \\
1 & 0 & 1 & 1 & 1 & 0 & 0 & 1 \\
1 & 0 & 1 & 0 & 1 & 1 & 1 & 0
\end{array}\right] $ & \small $ \left[\begin{array}{rrrrrrrr}
1 & 1 & 0 & 0 & 0 & 0 & 2 & 1 \\
1 & 1 & 0 & -1 & 0 & -2 & 0 & 0 \\
0 & 0 & 0 & -1 & 0 & -2 & -2 & -1 \\
0 & -1 & -1 & -3 & -1 & -5 & -1 & 0 \\
0 & 0 & 0 & -1 & 0 & -2 & -2 & -1 \\
0 & -2 & -2 & -5 & -2 & -8 & 0 & 1 \\
2 & 0 & -2 & -1 & -2 & 0 & 12 & 7 \\
1 & 0 & -1 & 0 & -1 & 1 & 7 & 4
\end{array}\right] $ 
\normalsize \TBstrut \\ \hline

8 & \verb+ Gv]mtW + & \small 
$ \left[\begin{array}{rrrrrrrr}
0 & 1 & 1 & 1 & 0 & 1 & 1 & 1 \\
1 & 0 & 0 & 1 & 1 & 0 & 1 & 0 \\
1 & 0 & 0 & 1 & 1 & 1 & 0 & 1 \\
1 & 1 & 1 & 0 & 1 & 0 & 1 & 0 \\
0 & 1 & 1 & 1 & 0 & 1 & 1 & 1 \\
1 & 0 & 1 & 0 & 1 & 0 & 0 & 1 \\
1 & 1 & 0 & 1 & 1 & 0 & 0 & 0 \\
1 & 0 & 1 & 0 & 1 & 1 & 0 & 0
\end{array}\right] $ & \small $ \left[\begin{array}{rrrrrrrr}
0 & -5 & 11 & -6 & 0 & 7 & -5 & 7 \\
-5 & -1 & 0 & -1 & -5 & 0 & -1 & 0 \\
11 & 0 & 5 & -1 & 11 & 3 & 0 & 3 \\
-6 & -1 & -1 & 1 & -6 & 0 & -1 & 0 \\
0 & -5 & 11 & -6 & 0 & 7 & -5 & 7 \\
7 & 0 & 3 & 0 & 7 & 2 & 0 & 2 \\
-5 & -1 & 0 & -1 & -5 & 0 & -1 & 0 \\
7 & 0 & 3 & 0 & 7 & 2 & 0 & 2
\end{array}\right] $ 
\normalsize \TBstrut \\ \hline

9 & \verb+ G~YmtW + & \small 
$ \left[\begin{array}{rrrrrrrr}
0 & 1 & 1 & 1 & 0 & 1 & 1 & 1 \\
1 & 0 & 1 & 1 & 1 & 0 & 1 & 0 \\
1 & 1 & 0 & 1 & 1 & 1 & 0 & 1 \\
1 & 1 & 1 & 0 & 0 & 0 & 1 & 0 \\
0 & 1 & 1 & 0 & 0 & 1 & 1 & 1 \\
1 & 0 & 1 & 0 & 1 & 0 & 0 & 1 \\
1 & 1 & 0 & 1 & 1 & 0 & 0 & 0 \\
1 & 0 & 1 & 0 & 1 & 1 & 0 & 0
\end{array}\right] $ & \small $ \left[\begin{array}{rrrrrrrr}
-2 & -3 & -3 & -1 & 0 & 1 & 1 & 1 \\
-3 & -3 & -4 & -1 & -1 & 0 & 1 & 0 \\
-3 & -4 & -3 & -2 & -1 & 1 & 0 & 1 \\
-1 & -1 & -2 & 0 & 0 & 0 & 1 & 0 \\
0 & -1 & -1 & 0 & 1 & 1 & 1 & 1 \\
1 & 0 & 1 & 0 & 1 & 1 & 0 & 1 \\
1 & 1 & 0 & 1 & 1 & 0 & 1 & 0 \\
1 & 0 & 1 & 0 & 1 & 1 & 0 & 1
\end{array}\right] $ 
\normalsize \TBstrut \\ \hline

10 & \verb+ GfymtW + & \small 
$ \left[\begin{array}{rrrrrrrr}
0 & 1 & 0 & 1 & 1 & 1 & 1 & 1 \\
1 & 0 & 0 & 1 & 1 & 0 & 1 & 0 \\
0 & 0 & 0 & 1 & 1 & 1 & 0 & 1 \\
1 & 1 & 1 & 0 & 0 & 0 & 1 & 0 \\
1 & 1 & 1 & 0 & 0 & 1 & 1 & 1 \\
1 & 0 & 1 & 0 & 1 & 0 & 0 & 1 \\
1 & 1 & 0 & 1 & 1 & 0 & 0 & 0 \\
1 & 0 & 1 & 0 & 1 & 1 & 0 & 0
\end{array}\right] $ & \small $ \left[\begin{array}{rrrrrrrr}
2 & 1 & 0 & 3 & 2 & 1 & 1 & 1 \\
1 & 1 & 0 & 1 & 1 & 0 & 1 & 0 \\
0 & 0 & -1 & -1 & -3 & -2 & 0 & -2 \\
3 & 1 & -1 & 4 & 0 & 0 & 1 & 0 \\
2 & 1 & -3 & 0 & -7 & -5 & 1 & -5 \\
1 & 0 & -2 & 0 & -5 & -3 & 0 & -3 \\
1 & 1 & 0 & 1 & 1 & 0 & 1 & 0 \\
1 & 0 & -2 & 0 & -5 & -3 & 0 & -3
\end{array}\right] $ 
\normalsize \TBstrut \\ \hline

11 & \verb+ GvY}tG + & \small 
$ \left[\begin{array}{rrrrrrrr}
0 & 1 & 1 & 1 & 0 & 1 & 1 & 1 \\
1 & 0 & 0 & 1 & 1 & 0 & 1 & 0 \\
1 & 0 & 0 & 1 & 1 & 1 & 0 & 1 \\
1 & 1 & 1 & 0 & 0 & 1 & 1 & 0 \\
0 & 1 & 1 & 0 & 0 & 1 & 1 & 0 \\
1 & 0 & 1 & 1 & 1 & 0 & 0 & 1 \\
1 & 1 & 0 & 1 & 1 & 0 & 0 & 0 \\
1 & 0 & 1 & 0 & 0 & 1 & 0 & 0
\end{array}\right] $ & \small $ \left[\begin{array}{rrrrrrrr}
4 & 3 & -1 & -17 & 0 & -1 & 3 & -3 \\
3 & 2 & 0 & -6 & 4 & 0 & 2 & 0 \\
-1 & 0 & 2 & 8 & 6 & 2 & 0 & 2 \\
-17 & -6 & 8 & 34 & 0 & 8 & -6 & 0 \\
0 & 4 & 6 & 0 & 17 & 6 & 4 & 0 \\
-1 & 0 & 2 & 8 & 6 & 2 & 0 & 2 \\
3 & 2 & 0 & -6 & 4 & 0 & 2 & 0 \\
-3 & 0 & 2 & 0 & 0 & 2 & 0 & -2
\end{array}\right] $ 
\normalsize \TBstrut \\ \hline

12 & \verb+ Gvw}tG + & \small 
$ \left[\begin{array}{rrrrrrrr}
0 & 1 & 1 & 1 & 1 & 0 & 1 & 1 \\
1 & 0 & 0 & 1 & 1 & 0 & 1 & 0 \\
1 & 0 & 0 & 1 & 1 & 1 & 0 & 1 \\
1 & 1 & 1 & 0 & 0 & 1 & 1 & 0 \\
1 & 1 & 1 & 0 & 0 & 1 & 1 & 0 \\
0 & 0 & 1 & 1 & 1 & 0 & 0 & 1 \\
1 & 1 & 0 & 1 & 1 & 0 & 0 & 0 \\
1 & 0 & 1 & 0 & 0 & 1 & 0 & 0
\end{array}\right] $ & \small $ \left[\begin{array}{rrrrrrrr}
4 & 3 & -2 & 1 & 2 & 0 & 3 & -1 \\
3 & 2 & 0 & 1 & 1 & 0 & 2 & 0 \\
-2 & 0 & 1 & 1 & -1 & 3 & 0 & -1 \\
1 & 1 & 1 & 1 & 0 & 1 & 1 & 0 \\
2 & 1 & -1 & 0 & 1 & -1 & 1 & 0 \\
0 & 0 & 3 & 1 & -1 & 1 & 0 & 1 \\
3 & 2 & 0 & 1 & 1 & 0 & 2 & 0 \\
-1 & 0 & -1 & 0 & 0 & 1 & 0 & -1
\end{array}\right] $ 
\normalsize \TBstrut \\ \hline

13 & \verb+ Gf{}tG + & \small 
$ \left[\begin{array}{rrrrrrrr}
0 & 1 & 0 & 1 & 1 & 0 & 1 & 1 \\
1 & 0 & 0 & 1 & 1 & 0 & 1 & 0 \\
0 & 0 & 0 & 1 & 1 & 1 & 0 & 1 \\
1 & 1 & 1 & 0 & 1 & 1 & 1 & 0 \\
1 & 1 & 1 & 1 & 0 & 1 & 1 & 0 \\
0 & 0 & 1 & 1 & 1 & 0 & 0 & 1 \\
1 & 1 & 0 & 1 & 1 & 0 & 0 & 0 \\
1 & 0 & 1 & 0 & 0 & 1 & 0 & 0
\end{array}\right] $ & \small $ \left[\begin{array}{rrrrrrrr}
-2 & 1 & 0 & 2 & 2 & 0 & 1 & -2 \\
1 & 1 & 0 & 1 & 1 & 0 & 1 & 0 \\
0 & 0 & 3 & 2 & 2 & 3 & 0 & -1 \\
2 & 1 & 2 & 2 & 2 & 2 & 1 & 0 \\
2 & 1 & 2 & 2 & 2 & 2 & 1 & 0 \\
0 & 0 & 3 & 2 & 2 & 3 & 0 & -1 \\
1 & 1 & 0 & 1 & 1 & 0 & 1 & 0 \\
-2 & 0 & -1 & 0 & 0 & -1 & 0 & -1
\end{array}\right] $ 
\normalsize \TBstrut \\ \hline

14 & \verb+ G~UlTW + & \small 
$ \left[\begin{array}{rrrrrrrr}
0 & 1 & 1 & 1 & 0 & 1 & 1 & 1 \\
1 & 0 & 1 & 1 & 1 & 0 & 0 & 0 \\
1 & 1 & 0 & 1 & 0 & 1 & 1 & 1 \\
1 & 1 & 1 & 0 & 1 & 0 & 0 & 0 \\
0 & 1 & 0 & 1 & 0 & 1 & 1 & 1 \\
1 & 0 & 1 & 0 & 1 & 0 & 0 & 1 \\
1 & 0 & 1 & 0 & 1 & 0 & 0 & 0 \\
1 & 0 & 1 & 0 & 1 & 1 & 0 & 0
\end{array}\right] $ & \small $ \left[\begin{array}{rrrrrrrr}
1 & 1 & 2 & 1 & 0 & 1 & -1 & 1 \\
1 & 1 & 3 & 1 & 1 & 0 & 0 & 0 \\
2 & 3 & 6 & 3 & 0 & 1 & -2 & 1 \\
1 & 1 & 3 & 1 & 1 & 0 & 0 & 0 \\
0 & 1 & 0 & 1 & -2 & 1 & -2 & 1 \\
1 & 0 & 1 & 0 & 1 & 1 & 0 & 1 \\
-1 & 0 & -2 & 0 & -2 & 0 & -1 & 0 \\
1 & 0 & 1 & 0 & 1 & 1 & 0 & 1
\end{array}\right] $ 
\normalsize \TBstrut \\ \hline

15 & \verb+ Gz[m~[ + & \small 
$ \left[\begin{array}{rrrrrrrr}
0 & 1 & 1 & 0 & 0 & 0 & 1 & 1 \\
1 & 0 & 1 & 1 & 1 & 0 & 1 & 1 \\
1 & 1 & 0 & 1 & 1 & 1 & 0 & 1 \\
0 & 1 & 1 & 0 & 1 & 0 & 1 & 0 \\
0 & 1 & 1 & 1 & 0 & 1 & 1 & 1 \\
0 & 0 & 1 & 0 & 1 & 0 & 1 & 1 \\
1 & 1 & 0 & 1 & 1 & 1 & 0 & 1 \\
1 & 1 & 1 & 0 & 1 & 1 & 1 & 0
\end{array}\right] $ & \small $ \left[\begin{array}{rrrrrrrr}
1 & 1 & 1 & 0 & 0 & 0 & 1 & 1 \\
1 & 2 & 2 & 1 & 1 & 0 & 2 & 1 \\
1 & 2 & 1 & 1 & -1 & -1 & 0 & -1 \\
0 & 1 & 1 & 1 & 1 & 0 & 1 & 0 \\
0 & 1 & -1 & 1 & -3 & -2 & -3 & -4 \\
0 & 0 & -1 & 0 & -2 & -1 & -2 & -2 \\
1 & 2 & 0 & 1 & -3 & -2 & -2 & -3 \\
1 & 1 & -1 & 0 & -4 & -2 & -3 & -3
\end{array}\right] $ 
\normalsize \TBstrut \\ \hline

16 & \verb+ Gv}lZ_ + & \small 
$ \left[\begin{array}{rrrrrrrr}
0 & 1 & 1 & 1 & 1 & 1 & 1 & 0 \\
1 & 0 & 0 & 1 & 1 & 0 & 0 & 1 \\
1 & 0 & 0 & 1 & 1 & 1 & 1 & 1 \\
1 & 1 & 1 & 0 & 1 & 0 & 0 & 1 \\
1 & 1 & 1 & 1 & 0 & 1 & 1 & 0 \\
1 & 0 & 1 & 0 & 1 & 0 & 1 & 0 \\
1 & 0 & 1 & 0 & 1 & 1 & 0 & 0 \\
0 & 1 & 1 & 1 & 0 & 0 & 0 & 0
\end{array}\right] $ & \small $ \left[\begin{array}{rrrrrrrr}
1 & 1 & 1 & 1 & 1 & 1 & 1 & 0 \\
1 & 1 & 0 & -1 & -1 & 0 & 0 & -1 \\
1 & 0 & 4 & 2 & 5 & 3 & 3 & 1 \\
1 & -1 & 2 & -7 & -1 & 0 & 0 & -3 \\
1 & -1 & 5 & -1 & 5 & 3 & 3 & 0 \\
1 & 0 & 3 & 0 & 3 & 2 & 2 & 0 \\
1 & 0 & 3 & 0 & 3 & 2 & 2 & 0 \\
0 & -1 & 1 & -3 & 0 & 0 & 0 & -1
\end{array}\right] $ 
\normalsize \TBstrut \\ \hline

17 & \verb+ GnulZ_ + & \small 
$ \left[\begin{array}{rrrrrrrr}
0 & 1 & 0 & 1 & 1 & 1 & 1 & 0 \\
1 & 0 & 1 & 1 & 1 & 0 & 0 & 1 \\
0 & 1 & 0 & 1 & 0 & 1 & 1 & 1 \\
1 & 1 & 1 & 0 & 1 & 0 & 0 & 1 \\
1 & 1 & 0 & 1 & 0 & 1 & 1 & 0 \\
1 & 0 & 1 & 0 & 1 & 0 & 1 & 0 \\
1 & 0 & 1 & 0 & 1 & 1 & 0 & 0 \\
0 & 1 & 1 & 1 & 0 & 0 & 0 & 0
\end{array}\right] $ & \small $ \left[\begin{array}{rrrrrrrr}
1 & 1 & 0 & 1 & 2 & 1 & 1 & 0 \\
1 & 1 & -3 & 1 & 2 & 0 & 0 & -1 \\
0 & -3 & -3 & -3 & 0 & 1 & 1 & -2 \\
1 & 1 & -3 & 1 & 2 & 0 & 0 & -1 \\
2 & 2 & 0 & 2 & 4 & 2 & 2 & 0 \\
1 & 0 & 1 & 0 & 2 & 2 & 2 & 0 \\
1 & 0 & 1 & 0 & 2 & 2 & 2 & 0 \\
0 & -1 & -2 & -1 & 0 & 0 & 0 & -1
\end{array}\right] $ 
\normalsize \TBstrut \\ \hline

18 & \verb+ Gvx|Ro + & \small 
$ \left[\begin{array}{rrrrrrrr}
0 & 1 & 1 & 1 & 1 & 0 & 1 & 0 \\
1 & 0 & 0 & 1 & 1 & 1 & 0 & 1 \\
1 & 0 & 0 & 1 & 1 & 1 & 1 & 1 \\
1 & 1 & 1 & 0 & 0 & 1 & 0 & 1 \\
1 & 1 & 1 & 0 & 0 & 1 & 1 & 1 \\
0 & 1 & 1 & 1 & 1 & 0 & 0 & 0 \\
1 & 0 & 1 & 0 & 1 & 0 & 0 & 0 \\
0 & 1 & 1 & 1 & 1 & 0 & 0 & 0
\end{array}\right] $ & \small $ \left[\begin{array}{rrrrrrrr}
1 & 2 & -1 & 4 & -3 & 0 & 1 & 0 \\
2 & 3 & 0 & 5 & -1 & 1 & 0 & 1 \\
-1 & 0 & -2 & 1 & -4 & -1 & 1 & -1 \\
4 & 5 & 1 & 8 & 0 & 2 & 0 & 2 \\
-3 & -1 & -4 & 0 & -7 & -2 & 1 & -2 \\
0 & 1 & -1 & 2 & -2 & 0 & 0 & 0 \\
1 & 0 & 1 & 0 & 1 & 0 & 1 & 0 \\
0 & 1 & -1 & 2 & -2 & 0 & 0 & 0
\end{array}\right] $ 
\normalsize \TBstrut \\ \hline

19 & \verb+ Gfx|Ro + & \small 
$ \left[\begin{array}{rrrrrrrr}
0 & 1 & 0 & 1 & 1 & 0 & 1 & 0 \\
1 & 0 & 0 & 1 & 1 & 1 & 0 & 1 \\
0 & 0 & 0 & 1 & 1 & 1 & 1 & 1 \\
1 & 1 & 1 & 0 & 0 & 1 & 0 & 1 \\
1 & 1 & 1 & 0 & 0 & 1 & 1 & 1 \\
0 & 1 & 1 & 1 & 1 & 0 & 0 & 0 \\
1 & 0 & 1 & 0 & 1 & 0 & 0 & 0 \\
0 & 1 & 1 & 1 & 1 & 0 & 0 & 0
\end{array}\right] $ & \small $ \left[\begin{array}{rrrrrrrr}
1 & 1 & 0 & 1 & 1 & 0 & 1 & 0 \\
1 & -4 & 0 & -3 & 1 & -1 & 0 & -1 \\
0 & 0 & 5 & -1 & 5 & 1 & 1 & 1 \\
1 & -3 & -1 & -2 & 0 & -1 & 0 & -1 \\
1 & 1 & 5 & 0 & 6 & 1 & 2 & 1 \\
0 & -1 & 1 & -1 & 1 & 0 & 0 & 0 \\
1 & 0 & 1 & 0 & 2 & 0 & 1 & 0 \\
0 & -1 & 1 & -1 & 1 & 0 & 0 & 0
\end{array}\right] $ 
\normalsize \TBstrut \\ \hline

20 & \verb+ G~^[\s + & \small 
$ \left[\begin{array}{rrrrrrrr}
0 & 1 & 1 & 1 & 0 & 1 & 1 & 1 \\
1 & 0 & 1 & 1 & 1 & 1 & 0 & 0 \\
1 & 1 & 0 & 1 & 1 & 0 & 0 & 1 \\
1 & 1 & 1 & 0 & 1 & 1 & 0 & 1 \\
0 & 1 & 1 & 1 & 0 & 1 & 1 & 1 \\
1 & 1 & 0 & 1 & 1 & 0 & 1 & 0 \\
1 & 0 & 0 & 0 & 1 & 1 & 0 & 1 \\
1 & 0 & 1 & 1 & 1 & 0 & 1 & 0
\end{array}\right] $ & \small $ \left[\begin{array}{rrrrrrrr}
1 & 1 & 4 & 3 & 0 & -1 & -1 & -1 \\
1 & 1 & 2 & 1 & 1 & -2 & 0 & 0 \\
4 & 2 & 8 & 6 & 4 & 0 & 0 & 2 \\
3 & 1 & 6 & 5 & 3 & 2 & 0 & 2 \\
0 & 1 & 4 & 3 & -2 & -2 & -2 & -3 \\
-1 & -2 & 0 & 2 & -2 & 7 & -1 & 0 \\
-1 & 0 & 0 & 0 & -2 & -1 & -1 & -2 \\
-1 & 0 & 2 & 2 & -3 & 0 & -2 & -3
\end{array}\right] $ 
\normalsize \TBstrut \\ \hline

21 & \verb+ G~X]\s + & \small 
$ \left[\begin{array}{rrrrrrrr}
0 & 1 & 1 & 1 & 0 & 0 & 1 & 1 \\
1 & 0 & 1 & 1 & 1 & 1 & 1 & 0 \\
1 & 1 & 0 & 1 & 1 & 0 & 0 & 1 \\
1 & 1 & 1 & 0 & 0 & 1 & 0 & 1 \\
0 & 1 & 1 & 0 & 0 & 1 & 1 & 1 \\
0 & 1 & 0 & 1 & 1 & 0 & 1 & 0 \\
1 & 1 & 0 & 0 & 1 & 1 & 0 & 1 \\
1 & 0 & 1 & 1 & 1 & 0 & 1 & 0
\end{array}\right] $ & \small $ \left[\begin{array}{rrrrrrrr}
1 & 1 & 2 & 1 & 0 & 0 & 1 & 5 \\
1 & 17 & -3 & 3 & 2 & -3 & 2 & 0 \\
2 & -3 & 5 & 1 & -1 & 0 & 0 & 11 \\
1 & 3 & 1 & 1 & 0 & -1 & 0 & 4 \\
0 & 2 & -1 & 0 & 0 & -1 & -1 & -1 \\
0 & -3 & 0 & -1 & -1 & -1 & -3 & 0 \\
1 & 2 & 0 & 0 & -1 & -3 & -4 & 3 \\
5 & 0 & 11 & 4 & -1 & 0 & 3 & 26
\end{array}\right] $ 
\normalsize \TBstrut \\ \hline

22 & \verb+ Grz[\{ + & \small 
$ \left[\begin{array}{rrrrrrrr}
0 & 1 & 1 & 0 & 1 & 1 & 1 & 1 \\
1 & 0 & 0 & 1 & 1 & 1 & 0 & 0 \\
1 & 0 & 0 & 1 & 1 & 0 & 0 & 1 \\
0 & 1 & 1 & 0 & 0 & 1 & 0 & 1 \\
1 & 1 & 1 & 0 & 0 & 1 & 1 & 1 \\
1 & 1 & 0 & 1 & 1 & 0 & 1 & 1 \\
1 & 0 & 0 & 0 & 1 & 1 & 0 & 1 \\
1 & 0 & 1 & 1 & 1 & 1 & 1 & 0
\end{array}\right] $ & \small $ \left[\begin{array}{rrrrrrrr}
4 & 1 & 2 & 0 & 2 & -4 & -1 & -1 \\
1 & 1 & 0 & 2 & 1 & 1 & 0 & 0 \\
2 & 0 & 1 & -1 & 2 & 0 & 0 & 1 \\
0 & 2 & -1 & 5 & 0 & 2 & 0 & -1 \\
2 & 1 & 2 & 0 & -4 & -14 & -3 & -7 \\
-4 & 1 & 0 & 2 & -14 & -24 & -5 & -15 \\
-1 & 0 & 0 & 0 & -3 & -5 & -1 & -3 \\
-1 & 0 & 1 & -1 & -7 & -15 & -3 & -8
\end{array}\right] $ 
\normalsize \TBstrut \\ \hline

23 & \verb+ G~xX{s + & \small 
$ \left[\begin{array}{rrrrrrrr}
0 & 1 & 1 & 1 & 1 & 0 & 0 & 1 \\
1 & 0 & 1 & 1 & 1 & 1 & 0 & 0 \\
1 & 1 & 0 & 1 & 1 & 0 & 1 & 0 \\
1 & 1 & 1 & 0 & 0 & 1 & 1 & 1 \\
1 & 1 & 1 & 0 & 0 & 1 & 1 & 1 \\
0 & 1 & 0 & 1 & 1 & 0 & 1 & 0 \\
0 & 0 & 1 & 1 & 1 & 1 & 0 & 1 \\
1 & 0 & 0 & 1 & 1 & 0 & 1 & 0
\end{array}\right] $ & \small $ \left[\begin{array}{rrrrrrrr}
-3 & 2 & 1 & -1 & -4 & 0 & 0 & -2 \\
2 & 5 & 2 & 3 & 5 & 1 & 0 & 0 \\
1 & 2 & 1 & 1 & 2 & 0 & 2 & 0 \\
-1 & 3 & 1 & 1 & 0 & 1 & -3 & -1 \\
-4 & 5 & 2 & 0 & -4 & 1 & -3 & -3 \\
0 & 1 & 0 & 1 & 1 & 1 & -4 & 0 \\
0 & 0 & 2 & -3 & -3 & -4 & 19 & -1 \\
-2 & 0 & 0 & -1 & -3 & 0 & -1 & -1
\end{array}\right] $ 
\normalsize \TBstrut \\ \hline

24 & \verb+ GvxX|s + & \small 
$ \left[\begin{array}{rrrrrrrr}
0 & 1 & 1 & 1 & 1 & 0 & 0 & 1 \\
1 & 0 & 0 & 1 & 1 & 1 & 0 & 0 \\
1 & 0 & 0 & 1 & 1 & 0 & 1 & 1 \\
1 & 1 & 1 & 0 & 0 & 1 & 1 & 1 \\
1 & 1 & 1 & 0 & 0 & 1 & 1 & 1 \\
0 & 1 & 0 & 1 & 1 & 0 & 1 & 0 \\
0 & 0 & 1 & 1 & 1 & 1 & 0 & 1 \\
1 & 0 & 1 & 1 & 1 & 0 & 1 & 0
\end{array}\right] $ & \small $ \left[\begin{array}{rrrrrrrr}
1 & 1 & 1 & 3 & 3 & 0 & 0 & 1 \\
1 & 0 & 0 & -2 & -2 & -1 & 0 & 0 \\
1 & 0 & 1 & 2 & 2 & 0 & 1 & 1 \\
3 & -2 & 2 & 0 & 0 & -1 & 4 & 2 \\
3 & -2 & 2 & 0 & 0 & -1 & 4 & 2 \\
0 & -1 & 0 & -1 & -1 & 0 & 1 & 0 \\
0 & 0 & 1 & 4 & 4 & 1 & 1 & 1 \\
1 & 0 & 1 & 2 & 2 & 0 & 1 & 1
\end{array}\right] $ 
\normalsize \TBstrut \\ \hline

25 & \verb+ GrxX[s + & \small 
$ \left[\begin{array}{rrrrrrrr}
0 & 1 & 1 & 0 & 1 & 0 & 0 & 1 \\
1 & 0 & 0 & 1 & 1 & 1 & 0 & 0 \\
1 & 0 & 0 & 1 & 1 & 0 & 1 & 0 \\
0 & 1 & 1 & 0 & 0 & 1 & 0 & 1 \\
1 & 1 & 1 & 0 & 0 & 1 & 1 & 1 \\
0 & 1 & 0 & 1 & 1 & 0 & 1 & 0 \\
0 & 0 & 1 & 0 & 1 & 1 & 0 & 1 \\
1 & 0 & 0 & 1 & 1 & 0 & 1 & 0
\end{array}\right] $ & \small $ \left[\begin{array}{rrrrrrrr}
1 & -4 & 3 & 0 & 1 & 0 & 0 & 3 \\
-4 & 25 & 0 & 3 & -4 & 15 & 0 & 0 \\
3 & 0 & 0 & 4 & -7 & 0 & -5 & 0 \\
0 & 3 & 4 & 1 & 0 & 5 & 0 & 4 \\
1 & -4 & -7 & 0 & -3 & -8 & -2 & -7 \\
0 & 15 & 0 & 5 & -8 & 9 & -4 & 0 \\
0 & 0 & -5 & 0 & -2 & -4 & -1 & -5 \\
3 & 0 & 0 & 4 & -7 & 0 & -5 & 0
\end{array}\right] $ 
\normalsize \TBstrut \\ \hline

26 & \verb+ Grx^]c + & \small 

$ \left[\begin{array}{rrrrrrrr}
0 & 1 & 1 & 0 & 1 & 0 & 1 & 1 \\
1 & 0 & 0 & 1 & 1 & 1 & 1 & 1 \\
1 & 0 & 0 & 1 & 1 & 0 & 1 & 0 \\
0 & 1 & 1 & 0 & 0 & 1 & 0 & 1 \\
1 & 1 & 1 & 0 & 0 & 1 & 1 & 0 \\
0 & 1 & 0 & 1 & 1 & 0 & 1 & 0 \\
1 & 1 & 1 & 0 & 1 & 1 & 0 & 1 \\
1 & 1 & 0 & 1 & 0 & 0 & 1 & 0
\end{array}\right] $ & \small $ \left[\begin{array}{rrrrrrrr}
0 & -1 & 1 & 0 & 2 & 0 & -2 & -1 \\
-1 & 0 & 0 & 1 & -1 & 1 & -3 & -1 \\
1 & 0 & 1 & 1 & 2 & 0 & 2 & 0 \\
0 & 1 & 1 & 4 & 0 & 2 & 0 & -1 \\
2 & -1 & 2 & 0 & 5 & -1 & 3 & 0 \\
0 & 1 & 0 & 2 & -1 & 1 & 1 & 0 \\
-2 & -3 & 2 & 0 & 3 & 1 & -11 & -4 \\
-1 & -1 & 0 & -1 & 0 & 0 & -4 & -1
\end{array}\right] $ 
\normalsize \TBstrut \\ \hline

27 & \verb+ GryW[k + & \small 
$ \left[\begin{array}{rrrrrrrr}
0 & 1 & 1 & 0 & 1 & 1 & 0 & 1 \\
1 & 0 & 0 & 1 & 1 & 0 & 0 & 0 \\
1 & 0 & 0 & 1 & 1 & 0 & 0 & 0 \\
0 & 1 & 1 & 0 & 0 & 1 & 0 & 1 \\
1 & 1 & 1 & 0 & 0 & 1 & 1 & 0 \\
1 & 0 & 0 & 1 & 1 & 0 & 1 & 1 \\
0 & 0 & 0 & 0 & 1 & 1 & 0 & 1 \\
1 & 0 & 0 & 1 & 0 & 1 & 1 & 0
\end{array}\right] $ & \small $ \left[\begin{array}{rrrrrrrr}
-1 & -1 & -1 & 0 & -1 & -1 & 0 & -1 \\
-1 & 0 & 0 & 1 & -1 & 0 & 0 & 0 \\
-1 & 0 & 0 & 1 & -1 & 0 & 0 & 0 \\
0 & 1 & 1 & 1 & 0 & 1 & 0 & 1 \\
-1 & -1 & -1 & 0 & 0 & 1 & 1 & 0 \\
-1 & 0 & 0 & 1 & 1 & 4 & 2 & 2 \\
0 & 0 & 0 & 0 & 1 & 2 & 1 & 1 \\
-1 & 0 & 0 & 1 & 0 & 2 & 1 & 1
\end{array}\right] $ 
\normalsize \TBstrut \\ \hline

28 & \verb+ GrW[[k + & \small 
$ \left[\begin{array}{rrrrrrrr}
0 & 1 & 1 & 0 & 0 & 0 & 1 & 1 \\
1 & 0 & 0 & 1 & 1 & 0 & 0 & 0 \\
1 & 0 & 0 & 1 & 1 & 0 & 0 & 0 \\
0 & 1 & 1 & 0 & 0 & 1 & 0 & 1 \\
0 & 1 & 1 & 0 & 0 & 1 & 1 & 0 \\
0 & 0 & 0 & 1 & 1 & 0 & 1 & 1 \\
1 & 0 & 0 & 0 & 1 & 1 & 0 & 1 \\
1 & 0 & 0 & 1 & 0 & 1 & 1 & 0
\end{array}\right] $ & \small $ \left[\begin{array}{rrrrrrrr}
1 & 3 & 3 & 0 & 0 & 0 & 5 & 4 \\
3 & 0 & 0 & 4 & -5 & 0 & 0 & 0 \\
3 & 0 & 0 & 4 & -5 & 0 & 0 & 0 \\
0 & 4 & 4 & 1 & 0 & 5 & 0 & -3 \\
0 & -5 & -5 & 0 & -1 & -4 & -3 & 0 \\
0 & 0 & 0 & 5 & -4 & 9 & -12 & -15 \\
5 & 0 & 0 & 0 & -3 & -12 & 16 & 20 \\
4 & 0 & 0 & -3 & 0 & -15 & 20 & 25
\end{array}\right] $ 
\normalsize \TBstrut \\ \hline

29 & \verb+ GvW[[k + & \small 
$ \left[\begin{array}{rrrrrrrr}
0 & 1 & 1 & 1 & 0 & 0 & 1 & 1 \\
1 & 0 & 0 & 1 & 1 & 0 & 0 & 0 \\
1 & 0 & 0 & 1 & 1 & 0 & 0 & 0 \\
1 & 1 & 1 & 0 & 0 & 1 & 0 & 1 \\
0 & 1 & 1 & 0 & 0 & 1 & 1 & 0 \\
0 & 0 & 0 & 1 & 1 & 0 & 1 & 1 \\
1 & 0 & 0 & 0 & 1 & 1 & 0 & 1 \\
1 & 0 & 0 & 1 & 0 & 1 & 1 & 0
\end{array}\right] $ & \small $ \left[\begin{array}{rrrrrrrr}
-3 & -1 & -1 & -1 & 0 & 0 & 1 & -4 \\
-1 & 0 & 0 & 1 & 1 & 0 & 0 & 0 \\
-1 & 0 & 0 & 1 & 1 & 0 & 0 & 0 \\
-1 & 1 & 1 & 6 & 0 & 1 & 0 & 5 \\
0 & 1 & 1 & 0 & 19 & -4 & -5 & 0 \\
0 & 0 & 0 & 1 & -4 & 1 & 1 & 1 \\
1 & 0 & 0 & 0 & -5 & 1 & 1 & 1 \\
-4 & 0 & 0 & 5 & 0 & 1 & 1 & 1
\end{array}\right] $ 
\normalsize \TBstrut \\ \hline

30 & \verb+ GrY[[k + & \small 
$ \left[\begin{array}{rrrrrrrr}
0 & 1 & 1 & 0 & 0 & 1 & 1 & 1 \\
1 & 0 & 0 & 1 & 1 & 0 & 0 & 0 \\
1 & 0 & 0 & 1 & 1 & 0 & 0 & 0 \\
0 & 1 & 1 & 0 & 0 & 1 & 0 & 1 \\
0 & 1 & 1 & 0 & 0 & 1 & 1 & 0 \\
1 & 0 & 0 & 1 & 1 & 0 & 1 & 1 \\
1 & 0 & 0 & 0 & 1 & 1 & 0 & 1 \\
1 & 0 & 0 & 1 & 0 & 1 & 1 & 0
\end{array}\right] $ & \small $ \left[\begin{array}{rrrrrrrr}
1 & 3 & 3 & 0 & 0 & 2 & 5 & 4 \\
3 & 0 & 0 & 4 & -5 & 0 & 0 & 0 \\
3 & 0 & 0 & 4 & -5 & 0 & 0 & 0 \\
0 & 4 & 4 & 1 & 0 & 1 & 0 & -3 \\
0 & -5 & -5 & 0 & -1 & -2 & -3 & 0 \\
2 & 0 & 0 & 1 & -2 & 1 & 4 & 5 \\
5 & 0 & 0 & 0 & -3 & 4 & 16 & 20 \\
4 & 0 & 0 & -3 & 0 & 5 & 20 & 25
\end{array}\right] $ 
\normalsize \TBstrut \\ \hline

31 & \verb+ Gv][[k + & \small 
$ \left[\begin{array}{rrrrrrrr}
0 & 1 & 1 & 1 & 0 & 1 & 1 & 1 \\
1 & 0 & 0 & 1 & 1 & 0 & 0 & 0 \\
1 & 0 & 0 & 1 & 1 & 0 & 0 & 0 \\
1 & 1 & 1 & 0 & 1 & 1 & 0 & 1 \\
0 & 1 & 1 & 1 & 0 & 1 & 1 & 0 \\
1 & 0 & 0 & 1 & 1 & 0 & 1 & 1 \\
1 & 0 & 0 & 0 & 1 & 1 & 0 & 1 \\
1 & 0 & 0 & 1 & 0 & 1 & 1 & 0
\end{array}\right] $ & \small $ \left[\begin{array}{rrrrrrrr}
13 & -3 & -3 & 4 & 0 & 4 & -5 & -2 \\
-3 & 0 & 0 & -1 & -1 & 0 & 0 & 0 \\
-3 & 0 & 0 & -1 & -1 & 0 & 0 & 0 \\
4 & -1 & -1 & 4 & 1 & 3 & 0 & 1 \\
0 & -1 & -1 & 1 & -1 & 2 & -1 & 0 \\
4 & 0 & 0 & 3 & 2 & 1 & 1 & 1 \\
-5 & 0 & 0 & 0 & -1 & 1 & 1 & 1 \\
-2 & 0 & 0 & 1 & 0 & 1 & 1 & 1
\end{array}\right] $ 
\normalsize \TBstrut \\ \hline

    % 23 & \verb+'G~xX{s'+ & \small $\begin{bmatrix*}[r]
    %                     0 & 1 & 1 & 1 & 1 & 0 & 0 & 1 \\
    %             1 & 0 & 1 & 1 & 1 & 1 & 0 & 0 \\
    %             1 & 1 & 0 & 1 & 1 & 0 & 1 & 0 \\
    %             1 & 1 & 1 & 0 & 0 & 1 & 1 & 1 \\
    %             1 & 1 & 1 & 0 & 0 & 1 & 1 & 1 \\
    %             0 & 1 & 0 & 1 & 1 & 0 & 1 & 0 \\
    %             0 & 0 & 1 & 1 & 1 & 1 & 0 & 1 \\
    %             1 & 0 & 0 & 1 & 1 & 0 & 1 & 0
    %                 \end{bmatrix*}$ & \small
    %                 $\begin{bmatrix*}[r]
    %                     -3 & 2 & 1 & -1 & -4 & 0 & 0 & -2 \\
    %                     2 & 5 & 2 & 3 & 5 & 1 & 0 & 0 \\
    %                     1 & 2 & 1 & 1 & 2 & 0 & 2 & 0 \\
    %                     -1 & 3 & 1 & 1 & 0 & 1 & -3 & -1 \\
    %                     -4 & 5 & 2 & 0 & -4 & 1 & -3 & -3 \\
    %                     0 & 1 & 0 & 1 & 1 & 1 & -4 & 0 \\
    %                     0 & 0 & 2 & -3 & -3 & -4 & 19 & -1 \\
    %                     -2 & 0 & 0 & -1 & -3 & 0 & -1 & -1
    %                 \end{bmatrix*}$ \normalsize \TBstrut \\ \hline
\end{longtable}
\end{center}

\bibliographystyle{plain}
\bibliography{MZ}

\end{document}